\documentclass[12pt,a4paper]{article}

\usepackage[left=1cm,right=1cm,top=2cm,bottom=2cm]{geometry}
\usepackage[dvips]{graphics}
\usepackage[dvips]{epsfig}
\usepackage{color}
\usepackage{hyperref}
\usepackage{tikz,authblk}

\newtheorem{lemma}{Lemma}
\newtheorem{theorem}[lemma]{Theorem}
\newtheorem{corollary}[lemma]{Corollary}
\newtheorem{proposition}[lemma]{Proposition}
\newtheorem{definition}{Definition}
\newtheorem{remark}{Remark}

\newcommand{\dimo}{\noindent \emph{Proof. }}
\newcommand{\qed}{\\ \rightline{$\Box$ \ \ \ \ \ \ \ \ \ \ \ \ \ \ \ }\\}
\newcommand{\e}{\varepsilon}
\newcommand{\G}{\Gamma}
\newcommand{\g}{\gamma}

\usepackage{latexsym}
\usepackage{amsfonts}
\usepackage{amssymb}
\usepackage{amsmath}
\usepackage[english]{babel}

\begin{document}

\title{Gem-induced trisections of compact PL $4$-manifolds}

 \renewcommand{\Authfont}{\scshape\small}
 \renewcommand{\Affilfont}{\itshape\small}
 \renewcommand{\Authand}{ and }

\author[1] {Maria Rita Casali}
\author[2] {Paola Cristofori}

\affil[1] {Department of Physics, Mathematics and Computer Science, University of Modena and Reggio Emilia \ \ \ \ \ \ \ \ \ \ \   Via Campi 213 B, I-41125 Modena (Italy), casali@unimore.it}

\affil[2] {Department of Physics, Mathematics and Computer Science, University of Modena and Reggio Emilia, Via Campi 213 B, I-41125 Modena (Italy), paola.cristofori@unimore.it}

\maketitle

\abstract{\noindent The idea of studying trisections of closed smooth $4$-manifolds via (singular) triangulations, endowed with a suitable vertex-labelling by three colors, is due to  Bell, Hass, Rubinstein and Tillmann, and has been applied by Spreer and Tillmann to standard simply-connected $4$-manifolds, via the so called {\it simple crystallizations}. 
In the present paper we propose generalizations of these ideas by taking into consideration a possible extension of trisections to compact PL $4$-manifolds with connected boundary, which is related to Birman's {\it special Heegaard sewing}, and by analyzing  {\it gem-induced trisections}, i.e.  trisections that can be induced not only by simple crystallizations, but also by any 5-colored graph encoding a PL  $4$-manifold with empty or connected boundary. 
This last notion gives rise to that of {\it G-trisection genus}, as an analogue, in this context, of the well-known trisection genus. 
We give conditions on a 5-colored graph ensuring one of its gem-induced trisections - if any - to realize  the G-trisection genus, and prove how to determine it directly from the graph. \\ 
As a consequence, we detect a class of closed simply-connected 4-manifolds, comprehending all standard ones, for which both G-trisection genus and trisection genus coincide with the second Betti number and also with half the value of the graph-defined PL invariant {\it regular genus}.\\
Moreover, the existence of gem-induced trisections and an estimation of the G-trisection genus via surgery description is obtained, for each compact PL 4-manifold admitting a handle decomposition lacking in 3-handles. 
}\endabstract



\section{Introduction}

The notion of {\it trisection}  of a  smooth, oriented closed $4$-manifold was introduced in 2016 by Gay and Kirby (\cite{Gay-Kirby}), by generalizing the classical idea of Heegaard splitting in dimension 3: via Morse functions, the $4$-manifold is decomposed into three $4$-dimensional handlebodies, with disjoint interiors and mutually intersecting in $3$-dimensional handlebodies, so that the intersection of all three ``pieces" is a closed orientable surface. The minimum genus of the intersecting surface is called the {\it trisection genus} of the $4$-manifold.  

In the last years, the notion of trisection has been extensively studied, also in its relationships with other manifold representation tools, such as handle-de\-com\-po\-si\-tions, Lefschetz fibrations,  
framed links and Kirby diagrams: see, for example,  \cite{Meier},  \cite{Castro-Ozbagci}, \cite{Bell-et-al}, \cite{SpreerTillmann}, \cite{Rubinstein-Tillmann} and \cite{Chu-Tillman}.  In particular, relying on the coincidence between DIFF and PL categories  
in dimension $4$, Bell, Hass, Rubinstein and Tillmann  faced  the study of unbalanced trisections  via (singular) triangulations, by making use of a vertex-labelling by three colors (a {\it tricoloring}, which  - under suitable conditions - actually encodes a trisection): see \cite{Bell-et-al} for details. Obviously, this approach brings all the advantages of a combinatorial description, and allows to algorithmically construct trisections  and estimate their genus, starting from singular triangulations.    

In \cite{SpreerTillmann} the same idea has been applied to the special case of {\it colored triangulations} associated to the so called {\it simple crystallizations} of simply-connected PL $4$-manifolds, thus succeeding in computing the trisection genus of all standard simply-connected $4$-manifolds (i.e. the connected sums of $\mathbb{CP}^2$ - possibly with reversed orientation -, $\mathbb S^2\times\mathbb S^2$, $\mathbb S^2\tilde\times\mathbb S^2$ and the $K3$-surface).  

The present paper performs a generalization of the ideas of  \cite{Bell-et-al} and \cite{SpreerTillmann}, along  
 two different directions: first, we take in consideration also compact PL $4$-manifolds with connected boundary, introducing a possible extension of the notion of trisection to the boundary case (restricted to compact $4$-manifolds whose associated singular manifold is simply-connected, and different from the one studied in \cite{Castro-Gay-Pinzon}); then, we prove that not only simple crystallizations, but all  edge-colored graphs encoding simply-connected $4$-manifolds give rise to a tricoloring, which always satisfies ``almost all" conditions to actually induce a trisection.
 
Note that, in order to study exotic structures, it is  important   to go beyond the class of simple crystallizations: in fact, as pointed out in  \cite{[Casali-Cristofori-Gagliardi JKTR 2015]}, there exist infinitely many simply-connected PL $4$-manifolds which are all TOP-homeomorphic  and do not admit simple crystallizations.   

Due to the features of generality of crystallization theory (i.e. the representation tool for compact PL manifolds of arbitrary dimension via edge-colored graphs) and related PL invariants  (mainly the {\it regular genus}, which extends to arbitrary dimension the classical genus of a surface and the Heegaard genus of closed 3-manifolds),  the approach performed in the present paper appears to be open to further developments, including a possible $n$-dimensional extension of the concept of trisection, based on a suitable 
partition of the vertices of the colored triangulation associated to a colored graph, for $n \ge 4$ (see \cite{Rubinstein-Tillmann} for a similar generalization to higher dimension). 
  
\bigskip 
\bigskip
    
The basic result connecting 4-dimensional crystallization theory (briefly reviewed in Section \ref{prelim}) with trisections is the following:

\begin{theorem}\label{thm.quasi-trisection}  Let $M^4$ be a compact PL $4$-manifold with empty (resp. connected)  boundary.   For each 5-colored graph 
$(\Gamma,\gamma)$ representing $M^4$  and for each cyclic permutation $\e=(\e_0, \dots, \e_4)$ of the color set   (resp. for each cyclic permutation $\e=(\e_0, \dots, \e_4)$ of the color set so that $\e_4$ is the only singular color  of $\G$), a triple $\mathcal  T(\Gamma, \e) =(H_{0},H_{1},H_{2})$ of $4$-dimensional submanifolds  of $M^4$ is obtained, so that: 
\begin{itemize}
 \item [(i)] $M^4 = H_{0}\cup H_{1}\cup H_{2}$ and the interiors of $H_{0},$ $H_{1},$ $H_{2}$ are pairwise disjoint;
 \item [(ii)]  $H_{1},H_{2}$ are $4$-dimensional handlebodies; $H_{0}$ is a $4$-disk (resp. is PL-homeomorphic to $\partial M^4 \times [0,1]); $
 \item [(iii)] $H_{01}=H_{0}\cap H_{1}$ and $H_{02}=H_{0}\cap H_{2}$ are $3$-dimensional handlebodies;
\item [(iv)]  $\Sigma (\mathcal T (\Gamma, \e)) = H_{0}\cap H_{1}\cap H_{2}$ is a closed connected surface (which is called \emph{central surface}), whose genus coincides with the regular genus, with respect to $\e,$ 
of the subgraph $\Gamma_{\hat{\e_4}}$, obtained from $\G$ by deleting all $\e_4$-colored edges.    
\end{itemize}
\end{theorem}

The triple $\mathcal T(\Gamma, \e)$  fits a weak notion of trisection, that we call {\it  quasi-trisection}: see Definition \ref{def_quasi-trisection}.
   
In case $H_{12}=H_{1}\cap H_{2}$ is a $3$-dimensional handlebody,   too, then $\mathcal T (\Gamma, \e)=(H_{0}, H_{1}, H_{2})$  is 
called a  {\it gem-induced trisection} of $M^4$ (see Definition \ref{def_B-trisection} and  Definition \ref{def_gem-induced trisection}   for details). 

\smallskip

As shown in Section \ref{quasi-trisections and B-trisections induced by colored graphs}, it is not difficult to check that, if $M^4$ admits a gem-induced trisection,  then the associated singular manifold $\widehat{M^4}$, obtained by capping off the connected boundary of $M^4$ by a cone over it,  is simply-connected.  Moreover, if $M^4$ is a closed 4-manifold, a gem-induced trisection 
is actually a trisection of $M^4$ (according to the definition by Gay and Kirby), with the property that one of the three involved 4-dimensional handlebodies is a 4-disk, while in the  connected  boundary case one of the 4-dimensional handlebodies is replaced by a collar of the boundary.  
Hence, the above notion may be considered as an extension of trisections also to compact PL $4$-manifolds with connected boundary (whose associated singular manifold is simply-connected), via an approach similar to the one suggested by Rubinstein and Tillmann in \cite[Subsection 5.2]{Rubinstein-Tillmann}, taking inspiration from an idea of Birman (\cite{Birman}). 

Moreover, note that  gem-induced trisections of closed $4$-manifolds are  $(g;k_1,k_2,0)$-trisections, which are studied in \cite{Meier_plus3}  and conjectured to exist for all closed simply-connected $4$-manifolds. 

 \bigskip 

By making use of the concept of gem-induced trisection, the following notion arises quite naturally: 

 \begin{definition}\label{def_GT-genus} {\em Let $M^4$ be a compact PL $4$-manifold with empty or connected  boundary, which admits at least a gem-induced trisection.  The  {\it  G-trisection genus} $g_{GT}(M^4)$ of $M^4$ is the minimum genus of the central surface of any gem-induced trisection of $M^4$: 
 $$   g_{GT}(M^4) \ = \  \min \{genus(\Sigma (\mathcal T (\Gamma, \e)))  \ /  \ \mathcal T (\Gamma, \e)  \ \ 
\text{is a gem-induced trisection of }  \ M^4\}.$$ 
}\end{definition}
  
\bigskip

In particular, in Section \ref{sec.Kirby-diagrams}, the existence of gem-induced trisections and an estimation of the G-trisection genus via surgery description is obtained, for each compact  PL 4-manifold admitting  a handle decomposition lacking in 3-handles. 

As it is well-known, any such $4$-manifold may be represented by a Kirby diagram, i.e. a link in $\mathbb S^3$ with some dotted components  and some framed components: dotted components describe  the attachment of 1-handles to the 4-disk having $\mathbb S^3$ as its boundary, framed components  describe  the attachment of 2-handles, while - only in case the 3-manifold obtained by Dehn surgery on the associated framed link  is $\mathbb S^3$ -  a  4-handle is trivially added so to  obtain the corresponding closed 4-manifold. 

\begin{theorem} \label{th_M4(L,c)} Let $M^4$ be a compact PL $4$-manifold with empty or connected boundary admitting a handle decomposition with no $3$-handles.  
\smallskip

\noindent 
Then, from each Kirby diagram $(L^{(*)},d)$ of $M^4$ a gem-induced trisection of $M^4$ can be algorithmically constructed, whose central surface has genus $s+1$, $s$ being the crossing number of the chosen diagram.

\noindent Furthermore, if  $(L^{(*)},d)$ has no dotted components (equivalently, if $1$-handles are missing, too), then a gem-induced trisection of $M^4$ can be obtained,  whose central surface has genus $m_\alpha$, $m_{\alpha}$ being the number of $\alpha$-colored regions in a chess-board coloration (by colors $\alpha$ and $\beta$) of the diagram.
\end{theorem}

Note that the class of compact PL 4-manifolds involved in the above Theorem is very large,  
and it is an open question whether it contains all closed simply-connected ones or not. In fact, in the closed case, the existence of 
a decomposition lacking in 3-handles (or, even, of a so called  {\it special handle decomposition},  lacking in 1-handles and 3-handles, according to \cite[Section 3.3]{[M]}) is related to Kirby problem n. 50, and is of particular interest with regard to exotic PL 4-manifolds: see, for example, \cite{[Ak1]} and  
\cite{[Ak2]}.

\medskip
Since in the closed case the G-trisection genus of $M^4$  gives an upper bound for the trisection genus $g_T(M^4)$, the following estimation of $g_T(M^4)$ follows:

\begin{corollary}\label{cor-framed-links} If a closed PL $4$-manifold $M^4$  admits a special handle decomposition, then: 
 $$g_T(M^4) \le  \min\{m_\alpha(L,c) \ /  \ (L,c)\ \text{is a framed link representing } M^4\}.$$
\end{corollary}

In particular, we obtain (Corollary \ref{linkfamilies}) that all $\mathbb D^2$-bundles over  $\mathbb S^2$ have G-trisection genus equal to $1$, thus proving that the G-trisection genus is not finite-to-one in the boundary case. On the other hand, it remains 
an open problem whether $g_{GT} $ is finite-to-one,  if a fixed (possibly empty) boundary is assumed. Note that this problem is related to Question 4.2 of \cite{Meier}, regarding the possible existence of infinitely many closed 4-manifolds admitting minimal $(g; 0)-$trisections, for some $g \ge 3$ (see also \cite{Lambert-Cole-Meier} and \cite{Meier-Zupan}).

Furthermore,  in Section \ref{sec.minimize_g_GT}, given  a compact 4-manifold $M^4$ (with empty or connected boundary and such that $\widehat{M^4}$ is simply-connected), we establish conditions on  a 5-colored graph representing  $M^4$ which ensure one of its gem-induced trisections - if any - to realize  the G-trisection genus of $M^4$, and prove how to determine it directly from the 
 graph: see Proposition \ref{minimal g_GT}  and Proposition \ref{minimality weak semi-simple}.  
Meanwhile, we detect a wide class of closed (resp. bounded) simply-connected $4$-manifolds -  containing all standard ones 
- for which the trisection genus (resp. the G-trisection genus), the second Betti number and half the regular genus coincide: see Proposition \ref{minimality weak semi-simple} and Corollary \ref{minimality closed weak simple}.

\section{Colored graphs encoding manifolds} \label{prelim}

The present section is devoted  to a brief  review of some basic notions of the so called {\it crystallization theory}, which is a representation tool for piecewise linear (PL) compact manifolds, without assumptions about dimension, connectedness, orientability or boundary properties  (see the ``classical" survey paper \cite{Ferri-Gagliardi-Grasselli}, or the more recent one \cite{Casali-Cristofori-Gagliardi Complutense 2015}, concerning the $4$-dimensional case). 

From now on, unless otherwise stated, all spaces and maps will be considered in the PL category, and all manifolds will be assumed to be compact, connected and orientable. 

\medskip

\begin{definition} \label{$n+1$-colored graph}
{\em An {\it $(n+1)$-colored graph}  ($n \ge 2$) is a pair $(\G,\g)$, where $\G=(V(\G), E(\G))$ is a multigraph (i.e. it can contain multiple edges, but no loops) 
which is regular of degree  $n+1$, and $\g$ is an {\it edge-coloration}, that is a map  $\g: E(\G) \rightarrow \Delta_n=\{0,\ldots, n\}$ which is injective on 
adjacent edges.}
\end{definition}

In the following, for sake of concision, when the coloration is clearly understood, we will denote the colored graph simply by $\G$. 

\smallskip

For every  $\{c_1, \dots, c_h\} \subseteq\Delta_n$ let $\G_{\{c_1, \dots, c_h\}}$  be the subgraph obtained from $\G$  by deleting all the edges that are not colored by the elements of $\{c_1, \dots, c_h\}$. 
In this setting, the complementary set of $\{c\}$ (resp. $\{c_1,\dots,c_h\}$)  in $\Delta_n$ will be denoted by $\hat c$ (resp. $\hat c_1\cdots\hat c_h$). 
The connected components of $\G_{\{c_1, \dots, c_h\}}$ are called {\it $\{c_1, \dots, c_h\}$-residues} or {\it $h$-residues} of $\G$; their number is denoted by $g_{\{c_1, \dots, c_h\}}$ (or, for short, by $g_{c_1,c_2}$, $g_{c_1,c_2,c_3}$ and $g_{\hat c}$ if $h=2,$ $h=3$ and $h = n$ respectively). 

 \medskip 

\noindent An $n$-dimensional pseudocomplex $K(\G)$ can be associated to an $(n+1)$-colored graph $\G$: 
\begin{itemize}
\item take an $n$-simplex for each vertex of $\G$ and label its vertices by the elements of $\Delta_n$;
\item if two vertices of $\G$ are $c$-adjacent ($c\in\Delta_n$), glue the corresponding $n$-simplices  along their $(n-1)$-dimensional faces opposite to the $c$-labeled vertices, so that equally labeled vertices are identified.
\end{itemize}

\smallskip

In general $|K(\G)|$ is an {\it $n$-pseudomanifold} and $\G$ is said to {\it represent} it. 

\medskip

Note that, by construction, $K(\G)$ is endowed with a vertex-labeling by $\Delta_n$ that is injective on any simplex. Moreover, $\G$ turns out to be the 1-skeleton of the dual complex of $K(\G)$.
The duality establishes a bijection between the $\{c_1, \dots, c_h\}$-residues of  $\G$  
and the $(n-h)$-simplices of $K(\G)$ whose vertices are labeled by   $\Delta_n - \{c_1, \dots, c_h\}$. 

Given a pseudocomplex $K$ and an $h$-simplex $\sigma^h$ of $K$, the {\it disjoint star} of $\sigma^h$ in $K$ is the pseudocomplex obtained by taking all $n$-simplices of $K$ having $\sigma^h$ as a face
and re-identifying only their faces that do not contain $\sigma^h.$ The {\it disjoint link}, $lkd(\sigma^h,K)$, of $\sigma^h$ in $K$ is the subcomplex of the disjoint star formed
by those simplices that do not intersect $\sigma^h.$

\noindent 

In particular, given an $(n+1)$-colored graph $\G$, each connected component of $\G_{\hat c}$ ($c\in\Delta_n$) is an $n$-colored graph representing the disjoint link of a $c$-labeled vertex of $K(\G)$, 
that is also (PL) homeomorphic to the link of this vertex in the first barycentric subdivision of $K.$

Therefore: 

\centerline{\it $|K(\G)|$ is a closed $n$-manifold iff, for each color $c\in\Delta_n$,}
\centerline{\it all $\hat c$-residues of $\G$ represent the $(n-1)$-sphere,}

while 

\centerline{\it $|K(\G)|$ is a singular\footnote{A {\it singular (PL) $n$-manifold} is a closed connected $n$-dimensional polyhedron admitting a simplicial triangulation where the links of vertices are closed connected $(n-1)$-manifolds. 
The notion extends also to polyhedra associated to colored graphs, provided that disjoint links of vertices are considered, instead of links.}   $n$-manifold iff, for each color $c\in\Delta_n$,} 
\centerline{\it  all $\hat c$-residues of $\G$ represent closed connected $(n-1)$-manifolds.}
\medskip

\begin{remark} \label{correspondence-sing-boundary} {\em If $N$ is a singular $n$-manifold, then a compact $n$-manifold $\check N$ is easily obtained by deleting small open neighbourhoods of its singular vertices.
Obviously $N=\check N$ iff $N$ is a closed manifold, otherwise $\check N$ has non-empty boundary (without spherical components).
Conversely, given a compact $n$-manifold $M$, a singular $n$-manifold $\widehat M$ can be constructed by capping off each component of $\partial M$ by a cone over it.

Note that, by restricting ourselves to the class of compact $n$-manifolds with no spherical boundary components,  the above correspondence is bijective and so singular $n$-manifolds and compact $n$-manifolds of this class can be associated  to each other in a well-defined way.

For this reason, throughout the present work, we will restrict our attention to compact manifolds without spherical boundary components. Obviously, in this wider context, closed $n$-manifolds are characterized by $M= \widehat M.$
}\end{remark}

In virtue of the bijection described in Remark \ref{correspondence-sing-boundary}, an $(n+1)$-colored graph $\G$ is said to {\it represent}
a compact $n$-manifold $M$ with no spherical boundary components (or, equivalently, to be a {\it gem} of $M$, where gem means {\it Graph Encoding Manifold})  if and only if  it represents the associated singular manifold $\widehat M$.

\medskip
The following theorem extends to the boundary case a well-known result - originally stated in \cite{Pezzana}  - founding the combinatorial representation theory for closed manifolds of arbitrary dimension via colored graphs. 

\begin{theorem}{\em (\cite{Casali-Cristofori-Grasselli})}\ \label{Theorem_gem}  
Any compact orientable (resp. non orientable) $n$-manifold with no spherical boundary components admits a bipartite (resp. non-bipartite) $(n+1)$-colored graph representing it.
\end{theorem}

If $\G$ is a gem of a compact $n$-manifold,  an $n$-residue of $\G$ will be called {\it ordinary} if it represents $\mathbb S^{n-1}$, {\it singular} otherwise. Similarly, a color $c$ will be called {\it singular} if at least one of the $\hat c$-residues of $\G$ is singular.

\bigskip

The  existence of a particular type of embedding of colored graphs into surfaces, is the key result in order to define the important notion of  regular genus. 
Since this paper concerns only orientable manifolds, we restrict the statement only to the bipartite case, although a similar result holds also for non-bipartite graphs. 

\begin{proposition}{\em (\cite{Gagliardi 1981})}\label{reg_emb}
Let $\G$ be a connected bipartite  
$(n+1)$-colored graph of order $2p$. Then for each cyclic permutation $\e = (\e_0,\ldots,\e_n)$ of $\Delta_n$, up to inverse, there exists a cellular embedding, called \emph{regular}, of $\G$  
into an orientable 
closed surface $F_{\e}(\G)$ whose regions are bounded by the images of the $\{\e_j,\e_{j+1}\}$-colored cycles, for each $j \in \mathbb Z_{n+1}$.
Moreover, the genus 
$\rho_{\e} (\G)$ of $F_{\e}(\G)$ satisfies

\begin{equation*}
2 - 2\rho_\e(\G)= \sum_{j\in \mathbb{Z}_{n+1}} g_{\e_j, \e_{j+1}} + (1-n)p.
\end{equation*}

\end{proposition}

\begin{definition} {\em The {\it regular genus} of  an $(n+1)$-colored graph $\G$ is defined as}
$$\rho(\G) = min\{\rho_\e(\G)\ |\ \e\ \text{cyclic permutation of \ } \Delta_n\}; $$
{\em the  {\it (generalized) regular genus} of a compact $n$-manifold $M$ is defined as}
$$\mathcal G (M) = min\{\rho(\G)\ |\ \G\ \text{gem of \ } M\}.$$
\end{definition}

\begin{remark} \label{rem_regular genus} {\em Note that the (generalized) regular genus is a PL invariant, extending to higher dimension the classical genus of a surface and the Heegaard genus of a $3$-manifold. It succeeds in 
characterizing spheres in arbitrary dimension (\cite{Ferri-Gagliardi Proc AMS 1982}), and a lot of classifying results via regular genus have been obtained, especially in dimension $4$ and $5$ (see \cite{Casali-Gagliardi ProcAMS},  \cite{Casali_Forum2003},  \cite{Casali-Cristofori-Gagliardi Complutense 2015}, \cite{generalized-genus}  and their references). Moreover, the regular genus is strictly related to the {\it G-degree}, a PL invariant arising within the theory of {\it Colored tensor models} in theoretical physics (\cite{Casali-Cristofori-Dartois-Grasselli}, \cite{Casali-Cristofori-Grasselli}, \cite{Casali-Grasselli 2019}).
}\end{remark}

\begin{definition} {\em An $(n+1)$-colored graph $\G$ representing a compact $n$-manifold with empty or connected boundary is said to be a {\it crystallization} of $M$ if,  for each color $c\in\Delta_n$, $\G_{\hat c}$ is connected. 
}\end{definition}

\noindent We recall that any compact $n$-manifold $M$ with empty or connected boundary admits a crystallization. In fact, given an arbitrary $(n+1)$-colored graph $\G$ representing $M$, it is always possible to perform a finite sequence of graph moves, called \textit{(proper) dipole moves}, which don't affect the represented manifold and at each step decreases the number of $\hat c$-residues of $\G$ for each $c\in\Delta_4$ (see \cite{Casali-Cristofori-Grasselli}).

\section{Trisections induced by colored graphs}
\label{quasi-trisections and B-trisections induced by colored graphs}

As already pointed out in Section  1, the fundamental result on which our approach to  trisections is based on consists in a weak notion of trisection ({\it quasi-trisection}), which arises from colored graphs representing compact $4$-manifolds. Then, it becomes quite natural to introduce  the notion of {\it B-trisection}, which actually extends the classical one of trisection and gives rise to that of gem-induced trisection. 
\medskip

We are now able to formalize the above notions.   

\begin{definition}\label{def_quasi-trisection}  {\em Let $M^4$ be a compact $4$-manifold with empty (resp. connected)  boundary.  A  {\it  quasi-trisection} 
of $M^4$  is a triple $\mathcal T =(H_{0},H_{1},H_{2})$ of $4$-dimensional submanifolds  of $M^4$,
whose interiors are pairwise disjoint, such that $M^4 = H_{0}\cup H_{1}\cup H_{2}$ and
\begin{itemize}
 \item [(i)]  $H_{1},H_{2}$ are $4$-dimensional handlebodies; $H_{0}$ is a $4$-disk \  (resp. is (PL) homeomorphic to $\partial M^4 \times [0,1]$);
 \item [(ii)] $H_{01}=H_{0}\cap H_{1}$ and $H_{02}=H_{0}\cap H_{2}$ are $3$-dimensional handlebodies;
\item [(iii)] 
$\Sigma (\mathcal T) = H_{0}\cap H_{1}\cap H_{2}$ is a closed connected surface (which is called {\it central surface}).  
\end{itemize}
}\end{definition}

\begin{definition}\label{def_B-trisection} {\em Let $M^4$ be a compact $4$-manifold with empty or connected  boundary. A  {\it  B-trisection} of $M^4$ is a quasi-trisection  $\mathcal T = (H_{0},H_{1},H_{2})$ of $M^4$  such that $H_{12}=H_{1}\cap H_{2}$ is a $3$-dimensional handlebody, too.}\end{definition}

 We use letter ``B" in B-trisections to point out that, if M is a closed (resp. bounded) 4-manifold, one of the 4-dimensional ``pieces" is a 4-Ball (resp. a collar on the Boundary).

\begin{remark} \label{rem-simply-connected} 
{\em Note that the central surface $\Sigma (\mathcal T)$ of a B-trisection $\mathcal T = (H_{0},H_{1},H_{2})$ of $M^4$ is a Heegaard surface for the 3-manifold $\partial H_i=\#_{k_i}  (\mathbb S^1 \times \mathbb S^2)$ ($k_i \ge 0$), for each $i\in \{1,2\}$, splitting it into the 3-dimensional handlebodies $H_{ij}$ and $H_{ik}$, with $\{j,k\}= \{0,1,2\}-\{i\}$. 
Moreover, in the closed (resp. boundary) case, $(H_{01},H_{02},\Sigma (\mathcal T))$ is a Heegaard splitting of $\partial H_0 = \mathbb S^3$ (resp. of $\partial M^4$, and more precisely of the boundary component of $\partial H_0$ intersecting $H_1\cup H_2$).
\\ \noindent
Hence, obviously, we have $ k_i \le genus(\Sigma (\mathcal T))$ for each $i\in \{1,2\}$, and $ \mathcal H (\partial M^4) \le genus(\Sigma (\mathcal T))$ in the boundary case (where $\mathcal H(M^3)$ denotes the classical Heegaard genus of the 3-manifold $M^3$).
\\  \noindent
Moreover, it is not difficult to check that, if $M^4$ admits a B-trisection, the associated singular manifold $\widehat{M^4}$ is simply-connected\footnote{Here, we obtain a generalization of the well-known fact that, if a closed $4$-manifold $M^4$ admits a trisection in which one of the ``pieces" is a $4$-disk, then $M^4$ is simply-connected.}: 
in fact, the existence of the above Heegaard splittings of $\partial H_1$ and $\partial H_2$, implies that there are surjections from $\pi_1(H_{12})$ to $\pi_1(\partial H_1)\cong\pi_1(H_1)$ 
and $\pi_1(\partial H_2)\cong\pi_1(H_2)$ respectively. Therefore, by applying Van-Kampen's theorem to the pair $(H_1,H_2)$, surjections are obtained $\pi_1(H_i)\twoheadrightarrow\pi_1(H_1\cup H_2)$ for each $i\in\{1,2\}.$
Now, if $H'_0$ denotes the cone over the boundary component of $H_0$ not intersecting $H_1\cup H_2$, then we have $(H_1\cup H_2)\cap H'_0 = H_{01}\cup H_{02} = \partial H'_0;$ 
since $\pi_1(\Sigma)\twoheadrightarrow\pi_1(H_{0i})\twoheadrightarrow\pi_1(\partial H_i)\cong\pi_1(H_i),\ \forall i\in\{1,2\}$, then there is a surjection 
$\pi_1(H_{01}\cup H_{02})\twoheadrightarrow\pi_1(H_1\cup H_2).$ Obviously, we have $H'_{0}\cup H_{1}\cup H_{2} = \widehat M^4$ and Van Kampen's theorem applied to the pair $(H_1\cup H_2,H'_0)$ implies $\pi_1(H'_0)\twoheadrightarrow\pi_1(\widehat M^4).$
Since $H'_0$ is contractible, $\widehat M^4$ is proved to be simply-connected.}
\end{remark} 
 
\begin{remark} \label{birman} {\em Note that, given a B-trisection $\mathcal T = (H_{0},H_{1},H_{2})$ of a 4-manifold $M^4$ with non-empty boundary, if both $H_1$ and $H_2$ are 4-balls, the Heegaard splitting of $\partial M^4$ induced by $\mathcal T$ defines a \textit{special Heegaard sewing} of this 3-manifold according to the definition in \cite{Birman}.
More precisely, in that paper Birman introduced a class of simply-connected 4-manifolds constructed by means of the partial identification of the boundaries of two 4-balls, in a way that, in our terminology, gives rise to a B-trisection
of the 4-manifold itself; consequently its boundary, which is a closed oriented 3-manifold, is said to be represented by a special Heegaard sewing.
Moreover, Birman proved that each closed orientable 3-manifold admits such a representation.
Rubinstein and Tillmann in \cite{Rubinstein-Tillmann} pointed out this construction  as a possible direction of generalization of trisections from closed to compact 4-manifolds with non-empty boundary.
As we will see in the next Proposition \ref{trisection}, this approach turns out to be quite natural starting from colored graphs.   
}\end{remark}

In order to prove Theorem \ref{thm.quasi-trisection}, a previous analysis on singular $4$-manifolds turns out to be useful. We will perform it by making use of the colored triangulations dual to colored graphs, and then by yielding from them a tricoloring,   
as in \cite{Bell-et-al} and \cite{SpreerTillmann}.

Let us denote by $G_s^{(4)}$ the set of $5$-colored graphs having only one $\hat 4$-residue and such that all $\hat i$-residues, with $i\in\Delta_3$, are ordinary.
Obviously, if $\G \in G_s^{(4)}, $  $|K(\G)|$ is a singular $4$-manifold with only one $4$-labelled, possibly singular, vertex, while all other vertices are not singular.

Observe that the class  $G_s^{(4)}$ properly contains (up to permutation of the color set) all crystallizations of compact $4$-manifolds with empty or connected boundary. Hence, any singular $4$-manifold $N^4$ with at most one singular vertex  can be represented by an element of $G_s^{(4)}$.

\medskip

From now on, when not otherwise stated, in the case of  a graph belonging to $G_s^{(4)}$ and while considering a cyclic  permutation $\e=(\e_0,\e_1,\e_2,\e_3,\e_4)$ of  $\Delta_4$, we will suppose $\e_4=4.$   On the other hand, given a gem of a closed simply-connected 4-manifold, the r\^ole of color 4 can be always played by any color $c$ such that $\G_{\hat c}$ is connected (in particular, by any color of $\Delta_4$ if $\G$ 
is a crystallization).

\begin{proposition}\label{trisection}
Given a singular $4$-manifold $N^4$ having at most one singular vertex, a $5$-colored graph $\G\in G_s^{(4)}$ representing $N^4$ and a cyclic permutation $\e$ of $\Delta_4$, there exists a triple $(H_{b},H_{r},H_{g})$  such that 
\begin{itemize}
\item [(i)] $H_{r},H_{g}$ are $4$-dimensional handlebodies; $H_{b}$ is either a $4$-ball, 
if $N^4$ is a closed $4$-manifold, or it is the cone over the disjoint link of the unique singular vertex of $N^4.$
\item [(ii)] $N^4 = H_{b}\cup H_{r}\cup H_{g}$ and the interiors of $H_{b},H_{r},H_{g}$ are pairwise disjoint.
\item [(iii)] $H_{rb}=H_{r}\cap H_{b}$ and $H_{gb}=H_{g}\cap H_{b}$ are $3$-dimensional handlebodies.
\item [(iv)] $H_{rg}=H_{r}\cap H_{g}$ collapses to a $2$-dimensional complex $\overline H(\G,\e).$
\item [(v)]  $H_{b}\cap H_{r}\cap H_{g}$ is a closed connected surface of genus $\rho_{\e_{\hat 4}}(\G_{\hat 4}),$ where $\e_{\hat 4} = (\e_0,\e_1,\e_2,\e_3).$ 
\end{itemize}

\end{proposition}

\dimo
Let us denote by $\sigma$ the standard $2$-simplex, by $v_b,v_r,v_g$ its vertices and by $\sigma^\prime$ its first barycentric subdivision. 

Following \cite{Bell-et-al} and \cite{SpreerTillmann}, let us consider the tricoloring of $K(\G)$ induced by the following partition of $\Delta_4:$ $b=\{\e_4=4\},\ g=\{\e_1,\e_3\},\ r=\{\e_0,\e_2\}.$ 

Along the present proof, for sake of simplicity, we suppose $\e = (0,1,2,3,4)$.

Let $\mu\ :\ K(\G)\to\sigma$ be the the simplicial map sending all equally colored vertices of $K(\G)$ to one vertex of the standard $2$-simplex; then $H_{b}$ (resp. $H_{r}$) (resp. $H_{g}$) is the preimage by $\mu$ of the star of $v_b$ (resp. $v_{r}$) (resp. $v_{g}$) in $\sigma^\prime.$ 
Therefore it is easy to see that $H_{b}$ is a regular neighbourhood of the (unique) $4$-colored vertex of $K(\G)$ and precisely it is the cone over its disjoint link. 
On the other hand $H_{r}$ (resp. $H_{g}$) is a regular neighbourhood of the $1$-dimensional subcomplex $K_{02}(\G)$ (resp. $K_{13}(\G)$) of $K(\G)$ generated by its 
$i$-colored vertices, with $i\in\{0,2\}$ (resp.  $i\in\{1,3\}$). 
Since $K(\G)$ is a pseudomanifold, $K_{02}(\G)$ (resp. $K_{13}(\G)$) is connected; therefore $H_{r}$ (resp. $H_{g}$) is a $4$-dimensional handlebody. 
 
This proves property (i). Property (ii) is straightforward.

\begin{figure}[t]
\centering
\scalebox{0.5}{\includegraphics{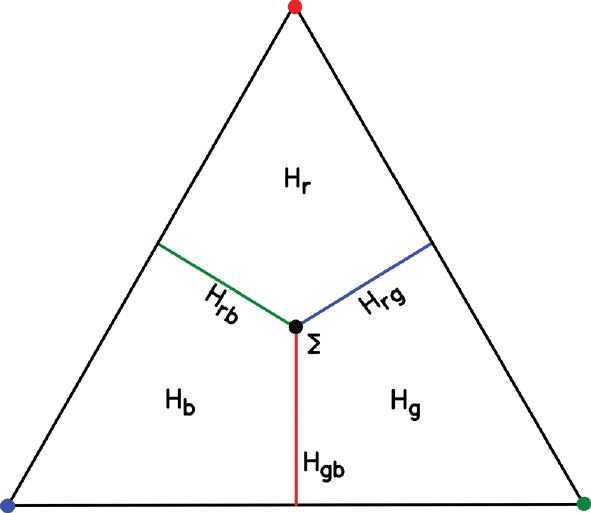}}
\caption{The 2-simplex $\sigma$}
\label{figsigma}
\end{figure}

\begin{figure}[t]
\centering
\scalebox{0.2}{\includegraphics{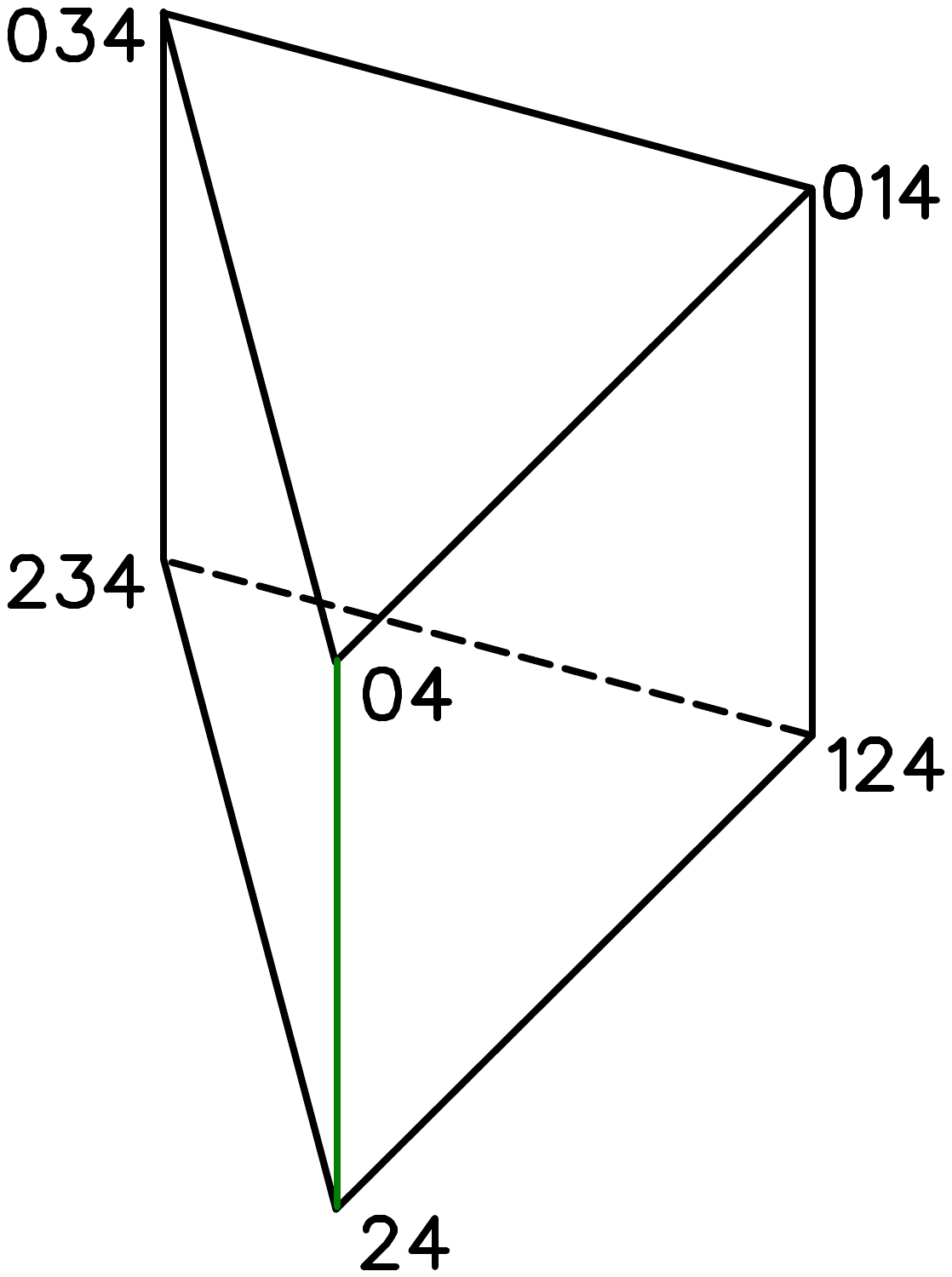}}
\caption{$H_{rb}$}
\label{figHrbprism}
\end{figure}

$H_{rb}$ is the preimage under $\mu$ of the edge of $\sigma^\prime$ with endpoints the barycenter of $\sigma$ and the barycenter of the edge opposite to $v_g$ (the ``green'' edge in Figure \ref{figsigma}); it can be obtained by taking one triangular prism for each $4$-simplex of $K(\G)$ - hence for each vertex of $\G$ - as in 
Figure \ref{figHrbprism}, where the vertices of the prism are barycenters of $1$- and $2$-simplexes of $K(\G)$. We denoted a barycenter with the colors of the vertices spanning it.

The prisms are glued to each other in the following way:

for each $i\in \{1,3\}$ and for each $i$-colored edge $e$ of $\G$, we identify the quadrangular (resp. triangular) faces of the two prisms associated to the endpoints of $e$ whose vertices do not contain the label $i$.

\begin{figure}[t]
\centering
\scalebox{0.26}{\includegraphics{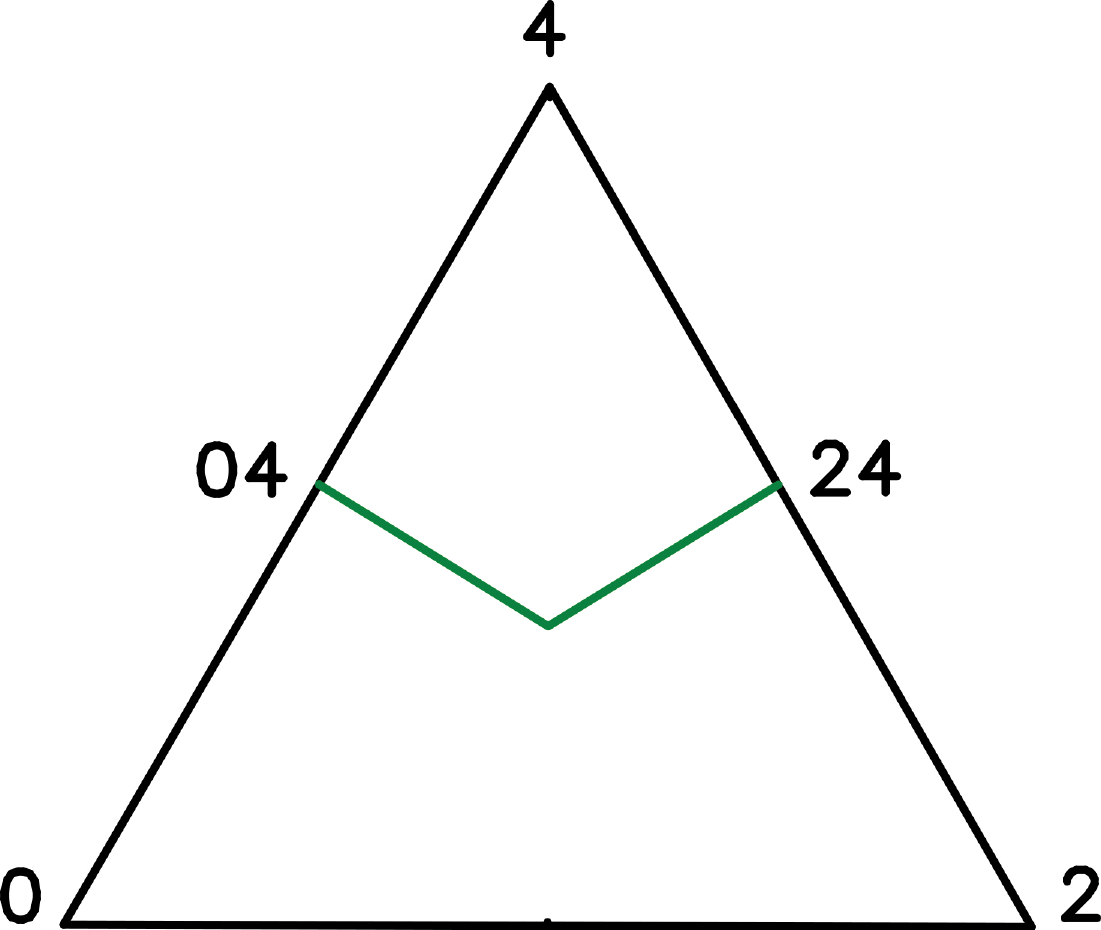}}
\caption{}
\label{figHrbarc}
\end{figure}

Note that the quadrangular face of each prism that is opposite to the ``green'' edge remains free.  
As a consequence, $H_{rb}$ collapses to the $1$-dimensional complex made by all the ``green'' edges, i.e. it collapses to a graph with $g_{1,2,3}+g_{0,1,3}$ vertices (i.e. the barycenters of the $\{2,4\}-$ and $\{0,4\}-$colored edges of $K(\G)$) and $g_{1,3}$ edges (each one is an arc inside a $\{0,2,4\}-$colored triangle of $K(\G)$ joining the barycenter of the triangle with the barycenters of the the $\{2,4\}-$ and $\{0,4\}-$colored edges; see Figure \ref{figHrbarc}).

Therefore $H_{rb}$ is a $3$-dimensional handlebody.    

\begin{figure}[t]
\centering
\scalebox{0.2}{\includegraphics{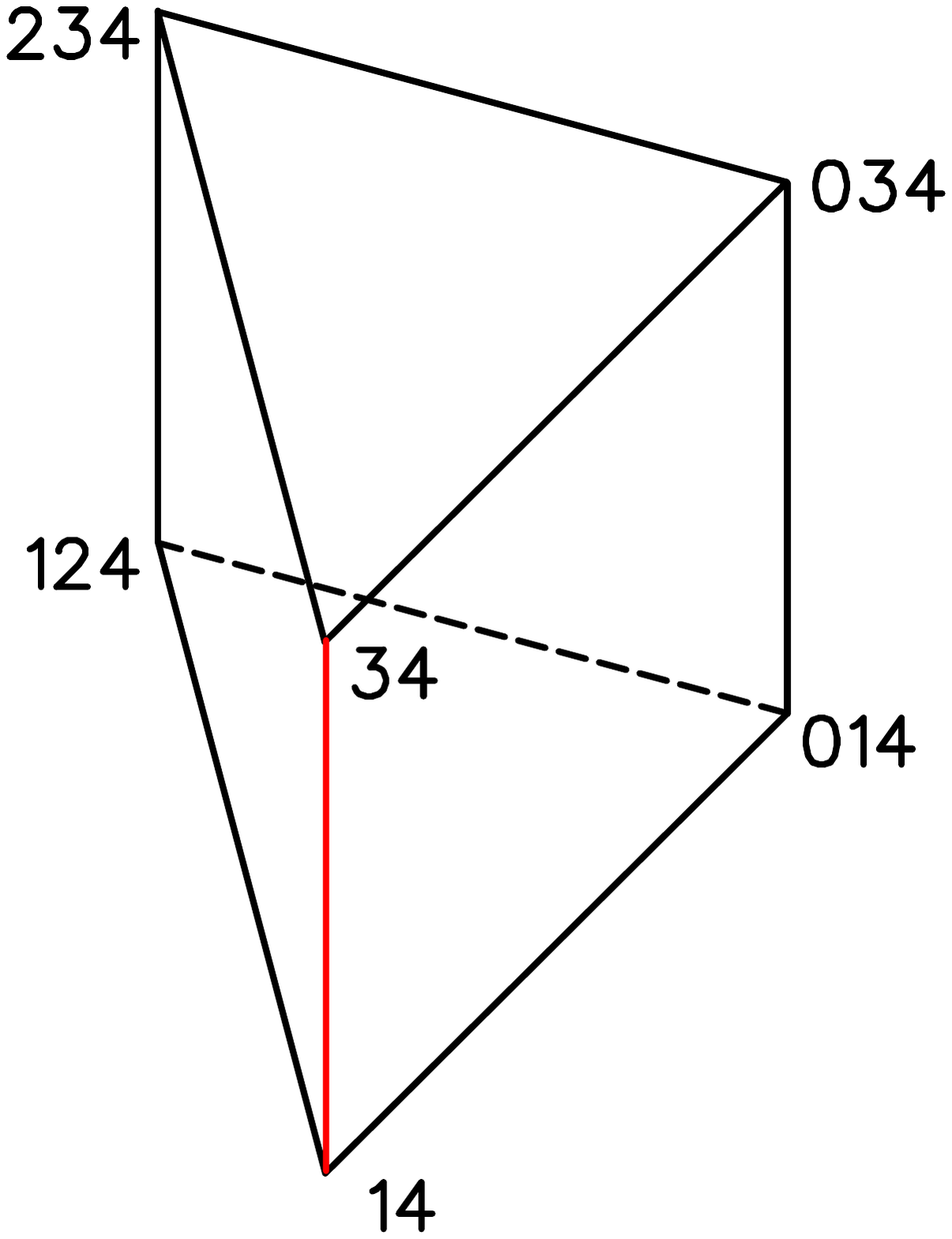}}
\caption{$H_{gb}$}
\label{figHgbprism}
\end{figure}

The proof that also $H_{gb}$ is a handlebody is analogous to the previous one with the exchange of colors $1,3$ with $0,2$ and of the vertex $v_g$ with $v_r$ (the ``red'' one in Figure \ref{figsigma}) of $\sigma$ (see Figure \ref{figHgbprism}).

\begin{figure}[t]
\centering
\scalebox{0.5}{\includegraphics{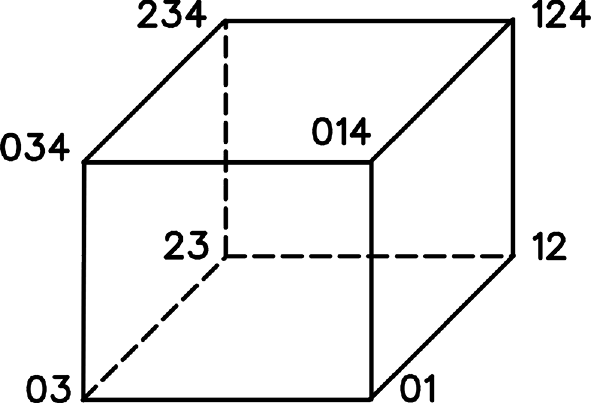}}
\caption{$H_{rg}$}
\label{figHrgcube}
\end{figure}

With regard to $H_{rg}$, this complex is the preimage under $\mu$ of the edge of $\sigma^\prime$ with endpoints the barycenter of $\sigma$ and the barycenter of the edge opposite to $v_b$ (the ``blue'' edge in Figure \ref{figsigma}); it can be obtained by taking a cube for each vertex of $\G$ as in Figure \ref{figHrgcube}, where, as above, the vertices are barycenters of $1$- and $2$-simplexes of $K(\G)$. 
The gluings of the cubes is performed according to the $i$-adjancencies of the vertices of $\G$, for each $i\in\Delta_4$, i.e. two cubes corresponding to the endpoints of an $i$-colored edge of $\G$ are glued along their faces whose vertices do not contain the color $i.$

\begin{figure}[t]
\centering
\scalebox{0.85}{\includegraphics{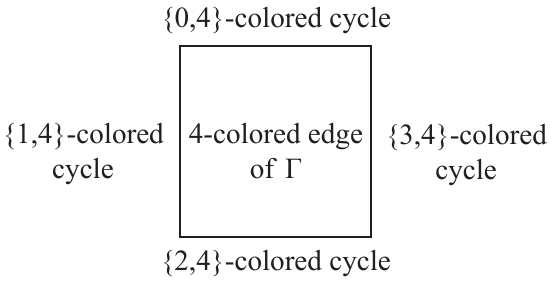}}
\caption{$\overline H(\G,\e)$}
\label{figHsquare}
\end{figure}

As a consequence, the face of each cube that misses no color (the top face in Figure \ref{figHrgcube}) is a free face and $H_{rg}$  
collapses to a $2$-dimensional complex $\overline H(\G,\e)$ consisting of one square for each $4$-colored edge of $\G$. It is not difficult to see that the edges of each square are dual to $\{i,4\}$-colored cycles of $\G$ for each $i\in\Delta_3$   (see Figure \ref{figHsquare}), therefore two squares are glued along an edge iff their corresponding $4$-colored edges of $\G$ share a common $\{i,4\}$-colored cycle of $\G.$

With regard to statement (v), note that $\Sigma = H_{b}\cap H_{r}\cap H_{g} = H_{rb}\cap H_{gb}\cap H_{rg}$ is a closed connected surface formed by the free faces of the cubes of $H_{rg}$, therefore by $2p$ squares, $4p$ edges and $g_{0,1}+g_{1,2}+g_{2,3}+g_{0,3}$ vertices, where $2p = \#V(\G)$. Therefore, by recalling that the subgraph $\G_{\hat 4}$ is connected, we have (via Proposition \ref{reg_emb}):
$$\chi(\Sigma) = g_{0,1}+g_{1,2}+g_{2,3}+g_{0,3}-2p =    2 - 2 \rho_{\e_{\hat 4}}(\G_{\hat 4}).$$ \qed

\medskip

By means of  Proposition \ref{trisection}, the existence of a quasi-trisection for each compact $4$-manifold  (Theorem \ref{thm.quasi-trisection}) easily 
follows. 
 
\medskip

{\it Proof of Theorem \ref{thm.quasi-trisection}.}\ \  Let $M^4$ be a compact $4$-manifold with empty or connected boundary, $\G$ a $5$-colored graph representing it and  $\e=(\e_0,\e_1,\e_2,\e_3,\e_4)$ a cyclic permutation of $\Delta_4$, with the condition that - if $\partial M^4\neq\emptyset$ - color $\e_4$ is the (only) singular color of $\G$. 

First of all note that, if $\G_{\widehat{\e_4}}$ is not connected, a suitable sequence of  (proper) dipole moves yields a new graph $\G'$, still representing $M^4$, having only one $\widehat{\e_4}$-residue   (see \cite{Casali-Cristofori-Grasselli}). 
In this case it is not difficult to prove that $\rho_{\e_{\hat 4}}(\G'_{\widehat{\e_4}}) = \rho_{\e_{\hat 4}}(\G_{\widehat{\e_4}})$, where $\rho_{\e_{\hat 4}}(\G_{\widehat{\e_4}})$ is intended as the sum of the regular genera of the connected components of $\G_{\widehat{\e_4}}$  with respect to the permutation $(\e_0,\e_1,\e_2,\e_3)$ induced by $\e$ on $\Delta_4 - \{\e_4\}.$ 

Therefore, without loss of generality, up to a permutation of colors, we can assume $\G\in G_s^{(4)}$ and $\e=(\e_0,\e_1,\e_2,\e_3,\e_4=4)$.

If $(H_{b},H_{r},H_{g})$ is the triple constructed as in Proposition \ref{trisection} starting from $\G$ and $\e$,  it is easy to check - via the bijection between $M^4$ and $\widehat{M^4}$ - that $(H_{b},H_{r},H_{g})$ defines a quasi-trisection $\mathcal T(\G,\e) = (H_{0},H_{1},H_{2})$ of $M^4$ where $H_{1}=H_{r}$ and $H_2=H_{g}$, while $H_0=H_b$, if $\partial M^4=\emptyset$, or $H_0$ is a collar of $\partial M^4$ obtained by removing from $H_b$ a suitable neighbourhood of the singular vertex,  if $\partial M^4 \ne \emptyset$.
Moreover, in both cases, $H_{0}\cap H_{1}\cap H_{2}= H_{b}\cap H_{r}\cap H_{g}$ is a closed connected surface of genus $\rho_{\e_{\hat 4}}(\G_{\widehat{\e_4}}).$
\qed
\bigskip

The notion of gem-induced trisection (from which the {\it G-trisection genus} is defined: see Definition \ref{def_GT-genus}) is used to identify B-trisections  arising from colored graphs.    

\begin{definition} \label{def_gem-induced trisection} {\em Let $M^4$ be a compact $4$-manifold $M^4$ with empty or connected  boundary, and let  $\mathcal T (\Gamma, \e)$ be the quasi-trisection of $M^4$ associated to the 5-colored graph $\G\in G_s^{(4)}$  and to the cyclic permutation $\e$ of $\Delta_4.$ 
In case $\mathcal T (\Gamma, \e)$ is a B-trisection, it is called  a {\it gem-induced trisection} of $M^4.$
}\end{definition}

Obviously,  in the particular case of closed simply-connected $4$-manifolds, the genus of the central surface of any gem-induced trisection gives an upper bound for the value of the ``classical" trisection genus. Hence, Theorem \ref{thm.quasi-trisection} directly yields the following statement:

\begin{corollary}
 If $\G$ is a crystallization of a closed simply-connected $4$-manifold $M^4$, so that $\mathcal T (\G,\e)$  is a gem-induced trisection ($\e$ being a cyclic permutation of $\Delta_4$), then the trisection genus of $M^4$ is less or equal to $\rho_{\e_{\hat 4}}(\G_{\widehat{\e_4}}).$
\ \qed 
\end{corollary}

Trisections of closed standard simply-connected $4$-manifolds induced by simple crystallizations are obtained in \cite{SpreerTillmann}, while gem-induced trisections for $4$-manifolds admitting handle decompositions lacking in $3$-handles will be constructed in the next section. 

The existence of a trisection induced by a given gem $\G\in G_s^{(4)}$ of  a compact $4$-manifold with empty or connected boundary is not always ensured; however it is possible to establish a combinatorial condition, which, if satisfied by a $5$-colored graph  $\G\in G_s^{(4)}$, implies that the quasi-trisection $\mathcal T (\G,\e)$ is actually  a B-trisection (and hence a gem-induced trisection), for each cyclic permutation $\e$ of $\Delta_4.$ 
  
\begin{proposition} \label{CS gem-induced trisections}
Let $\G\in G_s^{(4)}$ be a gem of  
a compact $4$-manifold $M^4$ with empty or connected boundary; with the notations of Proposition \ref{trisection}, for each cyclic permutation $\e$ of $\Delta_4$, the $2$-dimensional spine $\overline H(\G,\e)$ of $H_{rg}$  collapses to a graph if and only if  it is possible to establish an ordering $(e_1,\ldots,e_p)$ of the 4-colored edges of $\G$ such that for each $j\in\{1,\ldots,p\}$:
$$\begin{aligned}&\text{there exists } i\in\Delta_3 \text{ such that all 4-colored edges of  the } \{4,i\}\text{-colored cycle}\\
&\text{containing } e_j \text{ belong to the set }\{e_1,\ldots,e_j\}.\end{aligned}\hskip 80pt (*)$$
\noindent As a consequence, if the above condition is satisfied, then $\mathcal T(\Gamma, \e)$ is a gem-induced trisection of $M^4$, for each cyclic permutation $\e$ of $\Delta_4.$       
\end{proposition} 

\dimo  
First let us note that condition (*), for $j=1$, means that $e_1$ belongs,  for a suitable color $i$, to a $\{4,i\}$-colored cycle of $\G$ of length two;  hence, the square  in $\overline H(\G,\e)$ corresponding  to $e_1$ has a free edge from which it can be collapsed (see Figure \ref{figHsquare}). 

\noindent Moreover, it is not difficult to check that an ordering of the $4$-colored edges of $\G$ according to condition (*) exactly corresponds to a sequence of elementary collapses from  $\overline H(\G,\e)$ to a $1$-dimensional subcomplex. 
In fact, let us suppose condition (*) to be satisfied for $e_j$, $j\ge 2$, and for a color $i\in\Delta_3$, and let $C$ be the $\{4,i\}$-colored cycle containing $e_j$; then, after a sequence of elementary collapses (precisely those of the squares associated to the other $4$-colored edges of $C$, which belong to $\{e_1, \dots, e_{j-1}\}$), the square associated to $e_j$ has a free edge, namely the one corresponding to $C$, from which it can be collapsed.
Conversely, a sequence of elementary collapses from  $\overline H(\G,\e)$ to a graph naturally induces an ordering of the 4-colored edges of $\G$ (corresponding to the ordering of the collapsing squares of $\overline H(\G,\e)$) satisfying condition  (*). 
\noindent 
As regards the last statement, it is  sufficient to note that the collapse of all squares of $\overline H(\G,\e)$ proves that $H_{rg}$ is a 3-dimensional handlebody and, therefore,  that $\mathcal T(\Gamma, \e)$ is a gem-induced trisection of $M^4$, for each cyclic permutation $\e$ of $\Delta_4.$\qed

\begin{remark}\label{fund-group}{\em Within crystallization theory it is well-known the possibility of obtaining a presentation of the fundamental group of $\widehat{ M^4}$ whose generators are in bijection with the $4$-colored edges of a gem $\G\in G_s^{(4)}$ of $M^4$ and whose relators can be ``read'' on the $\{4,i\}$-colored cycles of $\G$ for all $i\in\Delta_3$ (see \cite{generalized-genus}). As a consequence,  condition (*) can be easily translated in terms of moves transforming this presentation into the trivial one.
\noindent More precisely, each move  is performed by taking a length one relation and deleting the involved generator from the generator set and from each relation containing it, too. The sequence of moves starts from the relation corresponding to the bicolored cycle of condition (*)  containing $e_1$,  which must have length two. } 
\end{remark}

\begin{remark} 
{\em As already pointed out above, in the case of crystallizations of closed simply-connected $4$-manifolds, quasi trisections can be constructed as in Proposition \ref{trisection} also by replacing color $4$ with any other color of $\Delta_4$; analogously condition (*) in Proposition \ref{CS gem-induced trisections} can be, more generally, formulated as depending on the choice of one color,  which is so associated to three distinct tricolorings of the dual triangulation.
\\
The implementation of this more general form of condition (*) into a C++ program, allows us to state that \underline{all} (rigid) crystallizations of closed $4$-manifolds up to 18 vertices induce trisections for all tricolorings associated to at least three colors (for the catalogue of such crystallizations see \cite{Casali-Cristofori ElecJComb 2015}). 
\\
The  program was also applied to two simple crystallizations obtained from the 16- and 17-vertices (possibly non PL-homeomorphic)  triangulations of the $K3$ surface constructed in \cite{Casella-Kuhnel}  and \cite{Spreer-Kuhnel} respectively;
both crystallizations turned out to induce trisections for all cyclic permutations of $\Delta_4.$
\\
As a particular case, we obtain again the result in  \cite[proof of Theorem 2]{SpreerTillmann} about trisections induced by simple\footnote{Simple crystallizations of a closed 
simply-connected $4$-manifold $M^4$ are characterized by having exactly $6\beta_2(M^4) + 2$ vertices (see \cite{[Casali-Cristofori-Gagliardi JKTR 2015]}).} crystallizations of all standard simply-connected $4$-manifolds.
}\end{remark}

\section{The case of $4$-manifolds not requiring 3-handles}\label{sec.Kirby-diagrams}

In this section we will take into account the case of $4$-manifolds, with empty or connected boundary, admitting a handle decomposition without 3-handles.  
It is well-known that any 4-manifold having such a decomposition  - apart from the genus $m$ handlebody $\mathbb Y^4_m$, with $m \ge0$   \footnote{A $5$-colored graph representing $\mathbb Y^4_m$, with $m>0$, can be found in \cite{generalized-genus} and we can consider it (together with the standard order two crystallization of $\mathbb D^4$) as associated to the Kirby diagram consisting only of $m \ge 0$ dotted components.} -   can be represented by a {\it  Kirby diagram}  $(L^{(m)},d)$, where $L$ is a link with $l$ components, $L_i$ being a dotted (resp. a framed) component $\forall i\in\{1,\ldots,m\}$, $0 \le m <l$,  (resp. 
$\forall i\in\{m+1,\ldots,l\}$), and $d=(d_1, \dots, d_{l-m})$, where $d_i \in \mathbb Z$, $\forall i \in \{1,\dots, l-m\}$, is the framing of the $(m+i)$-th (framed) component.
In fact, the $m$ dotted components represent the attachments of the 1-handles (see \cite{Akbulut-book} for this kind of notation), while the framed components of $(L^{(m)},d)$ describe the attachment of the $l-m$ 2-handles.  
Moreover, without loss of generality, we assume that the dotted components are unknotted, unlinked and, along them, overcrossings and undercrossings never alternate. 

As already recalled in the introduction, if $M^4(L^{(m)},d)$ denotes the $4$-manifold represented by $(L^{(m)},d)$,  its boundary  is the closed orientable 3-manifold $M^3(L,c)$ obtained from $\mathbb S^3$ by Dehn surgery along the associated framed link $(L,c)$, obtained by substituting each dotted component by a 0-framed one.
Moreover, in case $M^3(L,c)\cong \mathbb S^3$, we will still denote by $M^4(L^{(m)},d)$ the closed $4$-manifold obtained by adding a further $4$-handle.

\begin{figure}[!ht] 
\centerline{\scalebox{0.65}{\includegraphics{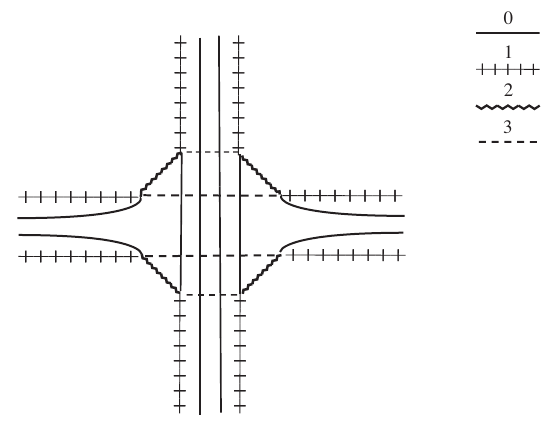}}}
\caption{\footnotesize 4-colored graph corresponding to a crossing}
\label{crossing}
\end{figure}

\begin{figure}[!ht] 
\centerline{\scalebox{0.50}{\includegraphics{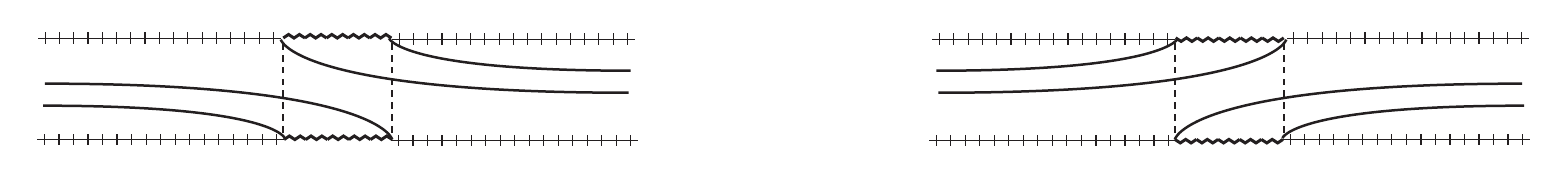}}}
\caption{\footnotesize 4-colored graphs corresponding to a positive curl (left) and a negative curl (right)}
\label{graph-curls}
\end{figure}

\begin{figure}[!h]
\centerline{\scalebox{0.4}{\includegraphics{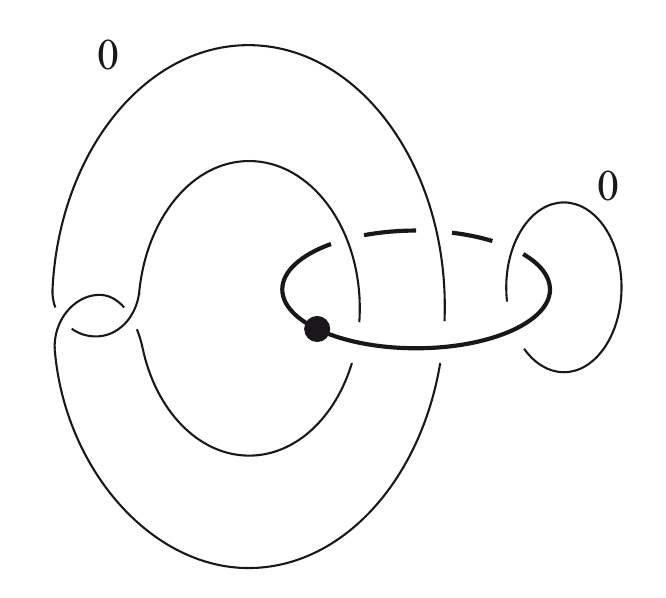}}}
\centerline{\scalebox{0.4}{\includegraphics{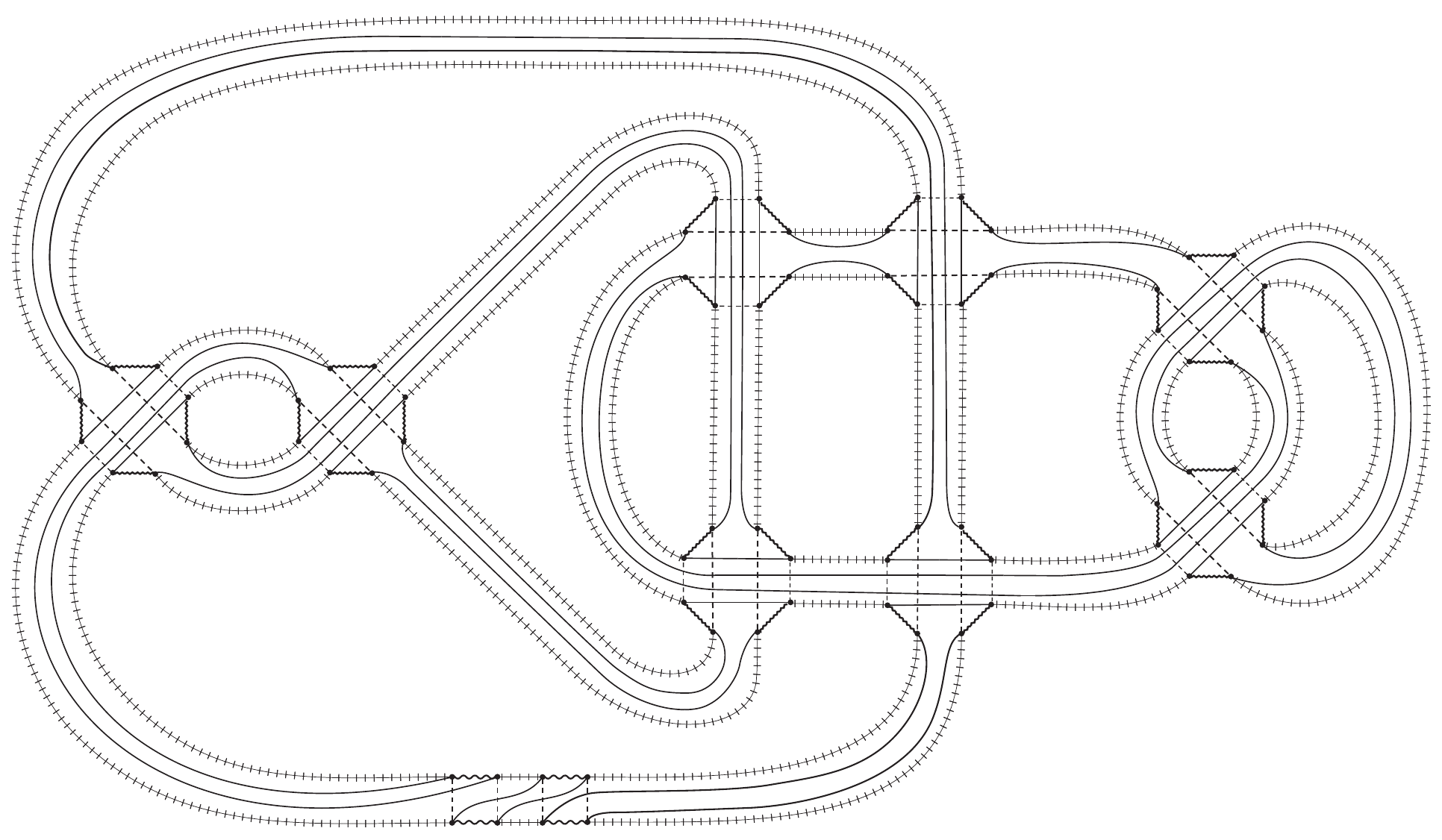}}}
\caption{\footnotesize A Kirby diagram and the associated 4-colored graph (representing $\mathbb S^2 \times \mathbb S^1$)}
\label{fig.fishtail3dim}
\end{figure}

In \cite{Casali-Cristofori Kirby-diagrams}  an algorithm is presented which enables to obtain a 5-colored graph representing $M^4(L^{(m)},d)$ starting from a planar diagram of $L$ with the possible addition of positive or negative curls in order to make the writhe of each component equal to its framing.
The first step of the procedure consists in considering a 4-colored graph  $\Lambda(L,c)$, representing the 3-manifold  $M^3(L,c)$, which can be directly ``drawn'' over $L$ as described in \cite{Casali JKTR2000}.  An example of this construction can be seen in Figure \ref{fig.fishtail3dim}, where the leftmost component of $L$ carries two additional negative curls.
As shown in the example, each crossing of $L$ gives rise to the graph shown in Figure \ref{crossing}, while each possible  curl gives rise to one of the graphs of Figure \ref{graph-curls}-left or Figure \ref{graph-curls}-right according to the curl being positive or negative.
Moreover, the hanging $0$- and $1$-colored edges of the above graphs are pasted together so that the $\{1,2\}$-colored cycles ``bound'' the regions of $L$ while $\{0,3\}$-colored cycles ``follow'' its components (see \cite{Casali JKTR2000} for details).

\begin{figure}[!h]
\centerline{\scalebox{0.45}{\includegraphics{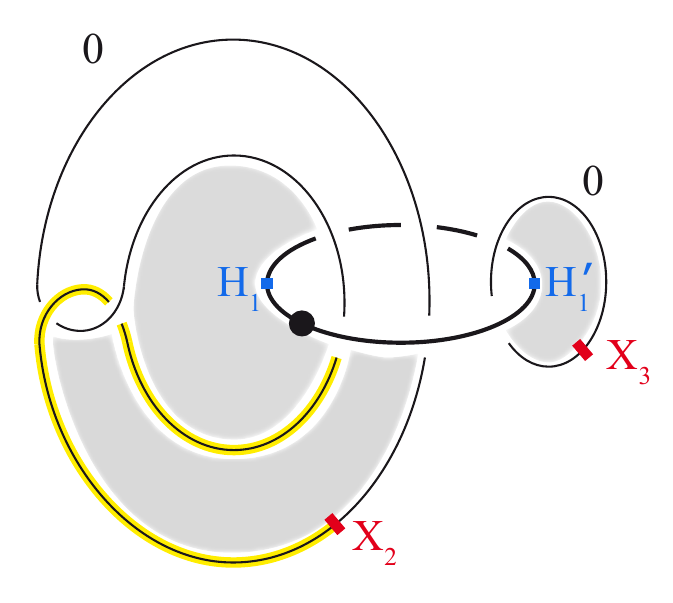}}}
\caption{\footnotesize the Kirby diagram of Fig. \ref{fig.fishtail3dim}, with points and highlighted segments of arcs.} 
\label{fig.fishtail_evidenziato}
\end{figure}

\begin{figure}[!h]
\centerline{\scalebox{0.45}{\includegraphics{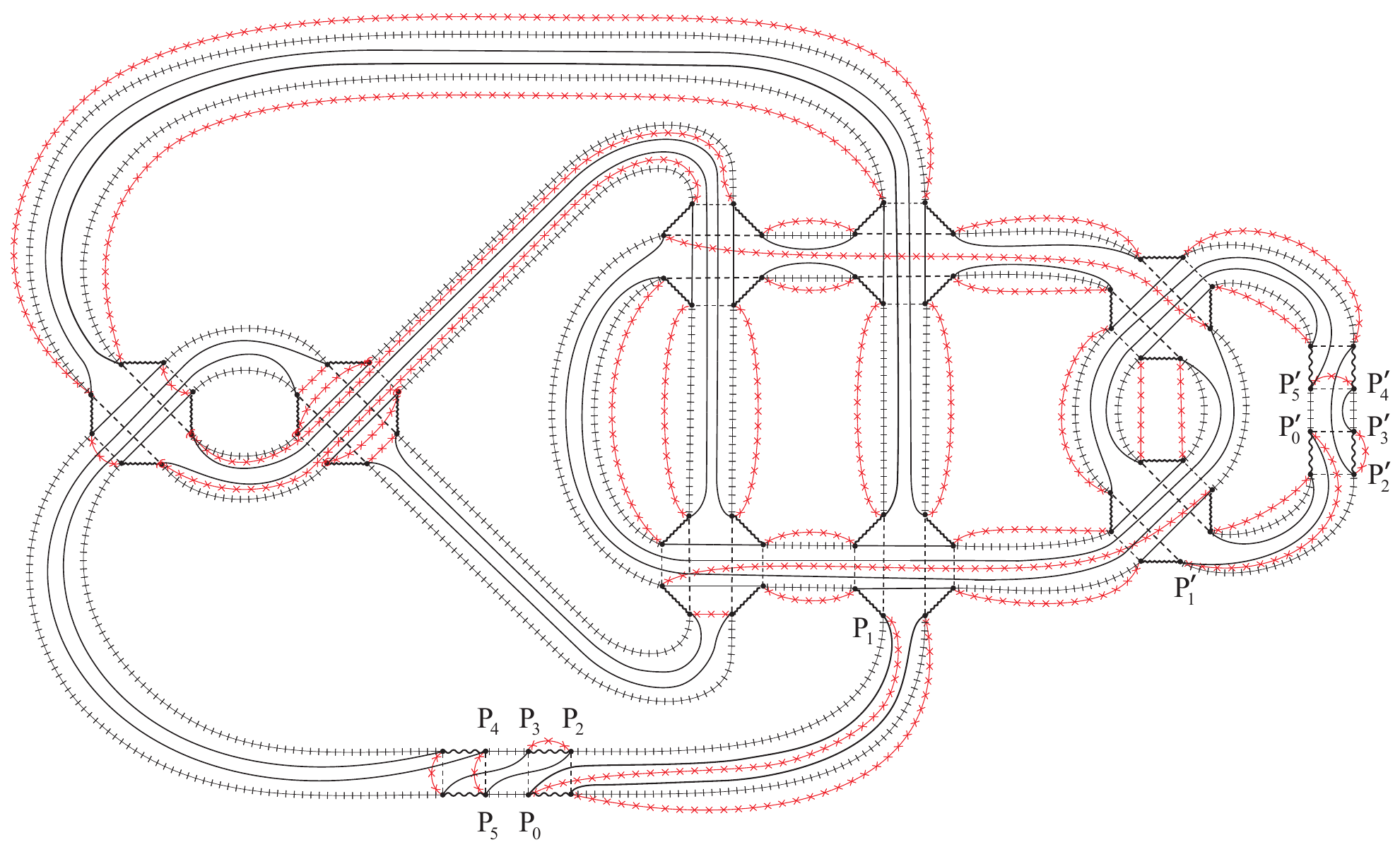}}}
\caption{\footnotesize The 5-colored graph associated to the Kirby diagram of Fig. \ref{fig.fishtail3dim} (representing $\mathbb S^2 \times \mathbb D^2$)}
 \label{fig.fishtail}
\end{figure}

Now, in order to construct a $5$-colored graph $\Gamma(L^{(m)},d)$ representing $M^4(L^{(m)},d)$, the following preliminary  choices and notations are needed:  
\begin{itemize}
\item[-] 
for each $i\in\{1,\ldots,m\},$ let $H_i$ and $H_i^{\prime}$ be  two points  of $L_i$ splitting it into one part containing only overcrossings and the other containing only undercrossings of $L_i$ (see Figure  \ref{fig.fishtail_evidenziato}).
\item[-] 
For each $j\in\{m+1,\ldots,l\},$ let us fix, possibly after adding to $L_j$ a trivial pair of opposite additional curls, a point $X_j$ of $L_j$ lying between a curl and an undercrossing\footnote{In the case of the trivial knot the role of the undercrossing can be played by another curl.}. Then a particular subgraph $Q_j$ of $\Lambda(L,c)$, called {\it quadricolor}, can be detected corresponding to the point $X_j$: see  Figure \ref{fig.quadricolor-singular}-left for the structure of a quadricolor, Figure \ref{fig.fishtail}  for an example (with a quadricolor for each framed component)  and \cite{Casali-Cristofori Kirby-diagrams} for details. 

\item[-] 
For each $j \in \{m+1, \dots l\}$ let us ``highlight" on $L_j$ - starting from $X_j$  and in the direction opposite to the undercrossing - a sequence $Y_j$ of consecutive  segments of arcs of $L- X_j$, so that  $H_i$ and $H_i^{\prime}$ belong to the boundary of the same region $\mathcal R_i$ of $L - \cup_{j=m+1}^l (X_j \cup Y_j)$, for each  $i\in\{1,\ldots,m\}$ (see Figure  \ref{fig.fishtail_evidenziato}, where  the yellow arcs form $Y_2$ and $Y_3$ is empty, while the shaded regions, together with the infinite one, give rise to $\mathcal R_1$).

Let us denote by $Y =Y_{m+1} \wedge \dots \wedge Y_l$ the sequence resulting from juxtaposition of all sequences of highlighted segments of arcs.

\item[-]
For each $i\in\{1,\ldots,m\}$,  let $\bar e_i$ (resp. $\bar e_i^{\prime}$) be the $1$-colored edge of $\Lambda(L,c)$ ``parallel" to the segment of arc of $L_i$ containing the point $H_i$ (resp. $H_i^{\prime}$)  ``on the side" of the regions of $L$  merging into $\mathcal R_i$   and let $v_i$ (resp.  $v_i^{\prime}$) be its end-point belonging to the subgraph corresponding to an undercrossing of the dotted component $L_i$.
\end{itemize}

The construction of $\Gamma(L^{(m)},d)$ consists in adding $4$-colored edges to  $\Lambda(L,c)$ according to  the following rules:

\begin{itemize}
 \item[(a)] for each $j\in\{m+1,\ldots, l\}$ and $\forall r \in \{0,1,2\}$, add a $4$-colored edge between the vertices $P_{2r}$ and $P_{2r+1}$ of the quadricolor $\mathcal Q_j$ (as shown in Figure \ref{fig.quadricolor-singular}); 
 
 \begin{figure}
\centerline{\scalebox{0.7}{\includegraphics{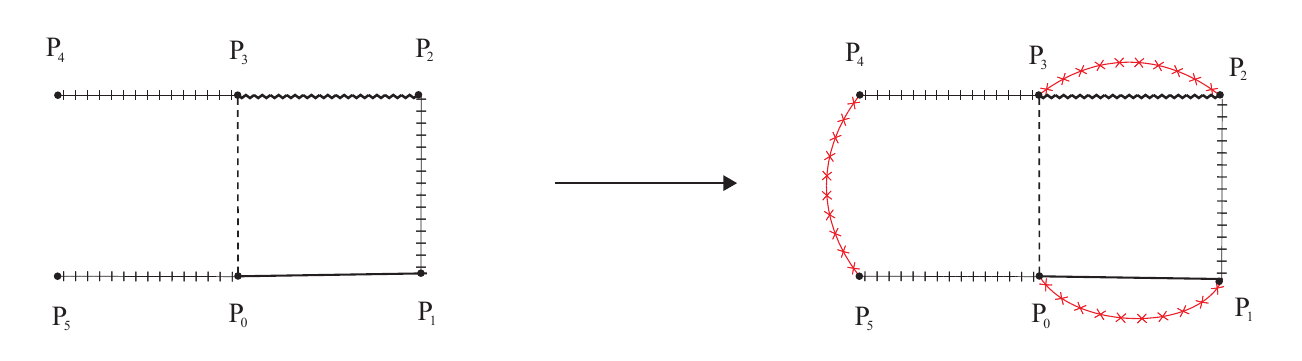}}}
\caption{\footnotesize step (a): adding edges to a quadricolor }
\label{fig.quadricolor-singular}
\end{figure}

\item[(b)] Follow the sequence $Y =Y_{m+1} \wedge \dots \wedge Y_l$, 
 starting, for each $j\in\{m+1,\ldots,l\}$ with $Y_j \ne \emptyset,$ from the segment of arc corresponding to the pair of $1$-colored edges adjacent to vertices $P_4$ and $P_5$ identified by the quadricolor $Q_j;$  at each step of the sequence, if  $f,f^{\prime}$ is the pair of $1$-colored edges which are ``parallel" to the considered highlighted segment of arc, then:

if no 4-colored edge has already been added to the end-points of $f$ and $f^\prime$, join, by $4$-colored edges, end-points of $f$ to end-points of $f^{\prime}$ belonging to different bipartition classes of $\Lambda(L,c);$ otherwise connect only the end-points of $f$ and $f^{\prime}$ having no already incident $4$-colored edge. 

Moreover, if two consecutive highlighted segments of arcs correspond to an undercrossing, whose overcrossing does not correspond to previous segments of arcs in $Y$, add $4$-colored edges so to double the pairs of $0$-colored edges within the subgraph corresponding to that crossing.

\item[(c)] For each $i\in\{1,\ldots, m\}$,  add a 4-colored edge, so to connect $v_i$ and ${v_i^{\prime}}$.

\item[(d)]  add $4$-colored edges between the remaining vertices of $\Lambda(L,c)$, joining those which belong to the same $\{1,4\}$-residue.

\end{itemize}

Figure \ref{fig.fishtail} shows the result of the above construction applied to the Kirby diagram in Figure \ref{fig.fishtail3dim}. 

\begin{proposition}{\em (\cite{Casali-Cristofori Kirby-diagrams})} For each Kirby diagram $(L^{(m)},d)$, the $5$-colored graph $\G(L^{(m)},d)$ represents the compact $4$-manifold $M^4(L^{(m)},d).$ 
\end{proposition}

\begin{figure}
\centerline{\scalebox{0.6}{\includegraphics{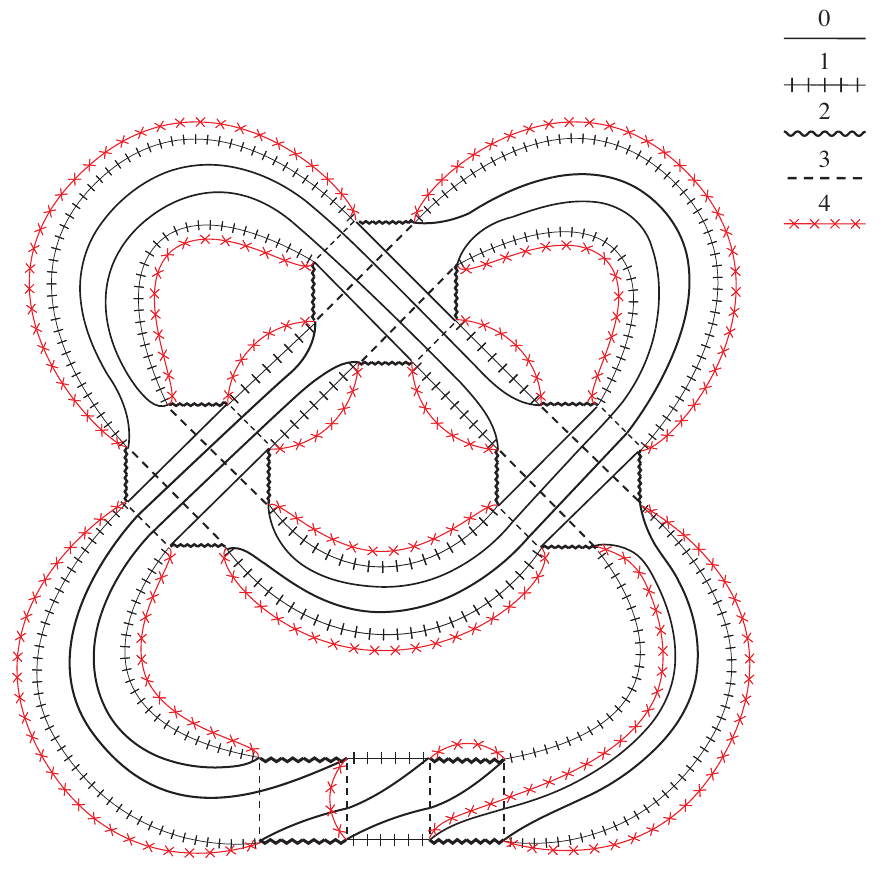}}}
\caption{\footnotesize The $5$-colored graph $\Gamma(L,c)$, for $c=+1$ and $L=$ trefoil} 
\label{fig.trefoil2}
\end{figure}

Now, we are able to prove Theorem \ref{th_M4(L,c)}, which ensures the existence of gem-induced trisections for all  
 $4$-manifolds with empty or connected boundary admitting a handle decomposition with no 3-handles.

\bigskip

{\it Proof of Theorem \ref{th_M4(L,c)}.}\ \  If $l=m  
$, it is easy to see that the $5$-colored graph in \cite{generalized-genus}, which represents the genus $m$  
handlebody,  belongs to $G_s^{(4)}$ and its unique ${\hat 4}$-residue has genus $m$ with respect to each cyclic permutation of $\Delta_3.$ Moreover, this $5$-colored graph trivially satisfies condition (*) of Proposition  \ref{CS gem-induced trisections} (actually all its $4$-colored edges belong to length two bicolored cycles), and therefore it induces a B-trisection for each cyclic permutation of $\Delta_4.$    

Let us now suppose $l-m>0$ and consider the $5$-colored graph  $\G(L^{(m)},d)$ constructed starting from a Kirby diagram $(L^{(m)},d)$ of $M^4$: note that 
$\G(L^{(m)},d)\in  G_s^{(4)},$ since  $4$ is its only possible singular color and $\G_{\hat 4}(L^{(m)},d)=\Lambda(L,c)$ is connected.    
In order to prove the first part of the statement, we will show that an ordering of the $4$-colored edges of  $\G(L^{(m)},d)$ can be arranged, so as to satisfy condition (*)  of Proposition \ref{CS gem-induced trisections}.

\medskip

First, we order the $4$-colored edges corresponding to each framed component in the following way.

For each $j\in\{m+1,\ldots,l\}$, let us consider the quadricolor $\mathcal Q_j$ and order the $4$-colored edges  $P_{2r}P_{2r+1}$,  added in step (a) of the construction, for increasing values of $r \in \{0,1,2\}$; note that  the first two edges belong to bicolored cycles of length two, while the $\{1,4\}$-colored cycle containing the third edge contains only the two previous ones. 

Let us then follow the sequence $Y_j$ and order the $4$-colored edges of step (b) as we meet them: condition (*) is still satisfied, by construction.

By going on along the component $L_j$, we then meet $4$-colored edges added in step (d): each of them belongs to a $\{1,4\}$-colored cycle whose other $4$-colored edges, if any, arise from step (b), so they already appear in the ordering sequence. 

By juxtaposing the obtained sequences of $4$-colored edges for increasing values of $j\in\{m+1,\ldots,l\}$, we get a sequence $\mathcal S$ of all $4$-colored edges of $\G(L^{(m)},d)$ corresponding to framed components. 

Now note that, for each $j\in\{1,\ldots,m\}$, any $\{1,4\}$-colored cycle involving $1$-colored edges which are  parallel to the (dotted)  component $L_j$, contains at most one $4$-colored edge not belonging to $\mathcal S$, with the exception of the $\{1,4\}$-colored cycle $C_j$ containing the vertices $v$ and $v'.$
Therefore we can add,  for each $j\in\{1,\ldots,m\}$, the $4$-colored edges of these cycles (different from $C_j$) to $\mathcal S$ in any order without breaking condition (*).

We point out that, after these additions,  the $4$-colored edges of $C_j$ not belonging to $\mathcal S$ are only two: one,  coming from step (c), shares the $\{3,4\}$-colored cycle (corresponding to the undercrossings of $L_j$) only with edges of $\mathcal S$, while the other, coming from step (d), shares the  $\{0,4\}$-colored cycle (corresponding to the overcrossings of $L_j$) only with edges of $\mathcal S$. 
As a consequence, by adding these last two edges to $\mathcal S$, for each $j\in\{1,\ldots,m\}$, we obtain the final ordering satisfying condition (*). 

This proves  - in virtue of Proposition \ref{CS gem-induced trisections} - that the quasi-trisection $\mathcal T(\G,\e)$  constructed starting from $\G(L^{(m)},d)$ is indeed a gem-induced trisection for each $\e\in\mathcal P_4.$ 

Moreover, if $\e = (1,0,2,3,4)$ is chosen, the central surface of $\mathcal T(\G,\e)$ has genus $\rho_{\e_{\hat 4}}(\G_{\hat 4}(L^{(m)},d))=\rho_{\e_{\hat 4}}(\Lambda(L,c))=s+1$,   $s$ being the number of crossings of $L$, as proved in \cite{Casali JKTR2000}.
\medskip

Let us now suppose $m=0$, i.e.  $(L^{(m)},d)=(L,c).$ 
Note that, in this case,  the  construction of $\G(L^{(m)},d)$ involves only steps (a) and (d); as a consequence, while the three $4$-colored edges involved in a quadricolor are as in Figure \ref{fig.quadricolor-singular}-right, all remaining $4$-colored edges double $1$-colored ones (see Figure \ref{fig.trefoil2}  for an example, where $L$ is the trefoil knot and $c=+1$).

Moreover, as proved in \cite{Casali-Cristofori Kirby-diagrams}, via a finite sequence of graph moves on $\Lambda(L,c)$ not affecting the quadricolor structures,  a new 4-colored graph  $\Omega(L,c)$ representing $M^3(L,c)$ can be considered whose regular genus with respect to the cyclic permutation $\e_{\hat 4}=(1,0,2,3)$ is $m_\alpha.$  
 Moreover,  a $5$-colored graph $\tilde \Omega(L,c)$, again representing $M^4(L,c)$, can be obtained by applying steps (a) and (d) to  $\Omega(L,c)$ instead of $\Lambda(L,c)$.  Figure \ref{fig.M4(K,c)} shows the result of the construction in the case of the trefoil knot with $c=+1$.

\begin{figure}[h] 
\centerline{\scalebox{0.5}{\includegraphics{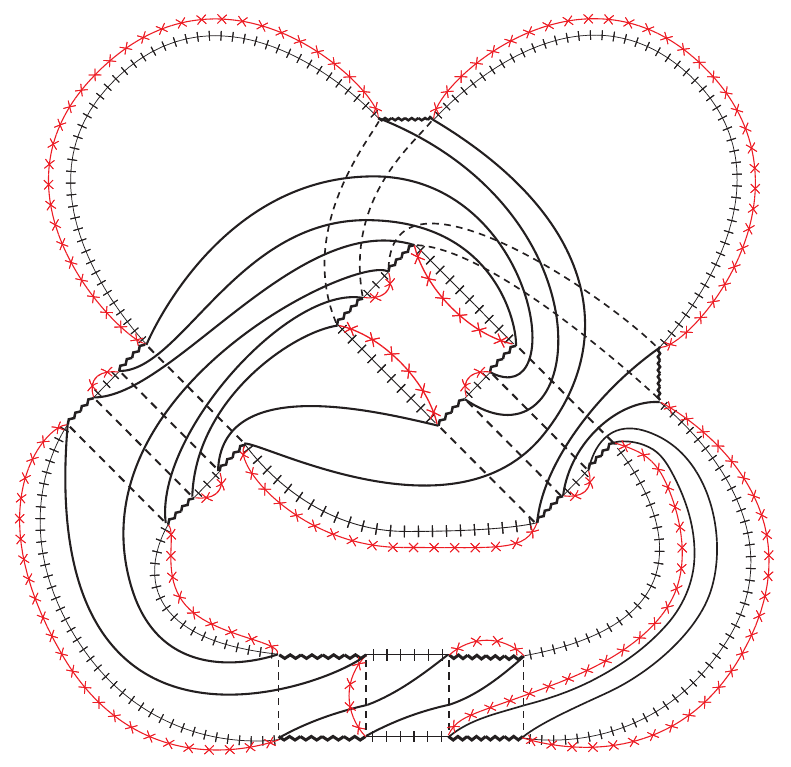}}}
\caption{\footnotesize  The $5$-colored graph $\tilde \Omega(L,c)$  representing  $M^4(L,c)$, for $c=+1$ and $L=$ trefoil}
\label{fig.M4(K,c)}
\end{figure}
 
Hence, it is easy to see that $\tilde \Omega(L,c)$  belongs to $G_s^{(4)}$ and satisfies condition (*);  therefore, by Proposition \ref{CS gem-induced trisections}, the quasi-trisection $\mathcal T(\tilde\Omega(L,c),\e)$ is a gem-induced trisection of $M^4(L,c).$
 
\smallskip
Finally, $\tilde \Omega_{\hat 4}(L,c)=\Omega(L,c)$ implies $\ genus(\Sigma(\mathcal T (\tilde\Omega(L,c),\e))= \rho_{\e_{\hat 4}}(\tilde \Omega_{\hat 4}(L,c)) = m_\alpha.$ \qed

\section{Results on the G-trisection genus}\label{sec-genus}

The results of Section \ref{quasi-trisections and B-trisections induced by colored graphs}
may be applied to yield information 
about the G-trisection genus (and hence, in the closed case, about the trisection genus). 

First of all,  the G-trisection genus is proved to satisfy the following properties:
 
\begin{proposition} \label{g_GT properties} \ \ 
\begin{itemize}  
\item[(i)]  $g_{GT}(M^4)=0 \ \ \Leftrightarrow  \ \ M^4 \cong \mathbb S^4$. 
\item[(ii)]  Let $M^4_1, M^4_2$ be two compact $4$-manifolds with empty or connected boundary which admit gem-induced trisections. 
If at least one between $M^4_1$ and $M^4_2$ is closed (resp. if both $M^4_1$ and $M^4_2$ have non-empty boundary), then $M^4_1 \# M^4_2$ (resp.  $M^4_1 \,^{\partial} \# M^4_2$) admits gem-induced trisections. \\ Moreover:   

$g_{GT}(M^4_1 \# M^4_2) \le g_{GT}(M^4_1) + g_{GT}(M^4_2)\quad$ 
(resp. $g_{GT}(M^4_1 \,^{\partial} \# M^4_2) \le g_{GT}(M^4_1) + g_{GT}(M^4_2)$) 
\item[(iii)] If $M^4$ has connected non-empty boundary,  then \ $g_{GT}(M^4) \ge \mathcal H(\partial M^4).$
\end{itemize}\end{proposition} 

\dimo First, it is easy to see that the standard order two crystallization of $\mathbb S^4$ induces a B-trisection of genus zero for each cyclic permutation of $\Delta_4.$ 
Conversely, if $g_{GT}(M^4)=0$, there exists  a $5$-colored graph $\G\in G_s^{(4)}$ representing $M^4$ and a cyclic permutation $\e$ of $\Delta_4$ so that $\mathcal T(\Gamma, \e)$ is a gem-induced trisection, whose central surface $\Sigma=\Sigma(\mathcal T (\Gamma,\e))$ is a 2-sphere.  
Since, by Theorem \ref{thm.quasi-trisection}, $\rho_{\e_{\hat 4}}(\G_{\widehat{\e_4}})=genus(\Sigma)=0$, then $\G_{\widehat{\e_4}}$ represents the 3-sphere; on the other hand $\G_{\widehat{\e_4}}$ is the unique possibly singular 4-residue of $\G$ and represents $\partial M^4$. Since, by Remark \ref{correspondence-sing-boundary}, spherical boundaries are not considered, then $M^4$ turns out to be closed. Moreover, by \cite[Proposition 7.3]{generalized-genus}, $M^4$ is either the $4$-sphere or a connected sum of copies of $\mathbb S^1\times\mathbb S^3$. The existence of a gem-induced trisection implies that $M^4$ is simply-connected, thus completing the proof of statement (i).

\smallskip
As concerns statement (ii), let us consider  two $5$-colored graphs $\G_1, \G_2 \in G_s^{(4)}$  and two cyclic permutations  $\e^{(1)},\ \e^{(2)}$ of $\Delta_4$ so that, for each $i \in \{1,2\}$, $\G_i$ represents $M^4_i$ and $\mathcal T(\Gamma_i, \e^{(i)})$  is a gem-induced trisection of $M^4_i$ with $g_{GT}(M^4_i) = genus (\Sigma(\mathcal T(\Gamma_i, \e^{(i)}))).$   
If necessary, perform a permutation on the color set of one of the two graphs, so to have $\e^{(1)} =\e^{(2)}=\bar \e.$

Note that, given a sequence of collapsing moves by which the submanifold $H_{rg}^{(i)}$ (made by the cubes in Figure \ref{figHrgcube}) of $|K(\G_i)|$ ($i\in\{1,2\}$) collapses to a graph, it is always possible to arrange it in order to perform first
the collapses of the 3-dimensional cubes. Let us choose now a vertex $v_1 \in V (\G_1)$ (resp. $v_2 \in V (\G_2)$) corresponding to a cube containing the 2-dimensional face of $H_{rg}^{(1)}$ that collapses last (resp. of $H_{rg}^{(2)}$ that collapses first) and let us perform the {\it graph connected sum} between $\G_1$ and $\G_2$ with respect to vertices $v_1$ and $v_2$ (i.e. the graph $\G_1 \#_{v_1,v_2} \G_2$ 
obtained by deleting $v_i$ in $\G_i$, $\forall i\in \{1,2\}$, and by welding the hanging edges of the same color: see  \cite[Section 7]{Grasselli-Mulazzani} for details about the construction).
$\G_1 \#_{v_1,v_2} \G_2$ is known to represent $M^4_1 \# M^4_2$  in case at least one between $M^4_1$ and $M^4_2$ is closed, and  $M^4_1 \,^{\partial} \# M^4_2$ in case both $M^4_1$ and $M^4_2$ have connected non-empty 
boundary.   

It is not difficult to see that the quasi-trisection induced by the graph $\G_1 \#_{v_1,v_2} \G_2$ and the cyclic permutation $\bar\e$ is indeed a B-trisection, since the collapsing sequences of $H_{rg}^{(1)}$ and $H_{rg}^{(2)}$ can be used to get a collapsing sequence of $H_{rg}^{(1)}\#H_{rg}^{(2)}$. Moreover,  the central surface of this gem-induced trisection is the connected sum of  $\Sigma(\mathcal T (\Gamma_1,\bar \e))$ and $\Sigma(\mathcal T (\Gamma_2,\bar \e)).$ 

In particular, note that, in case the collapsing sequence of  $H_{rg}^{(i)}$ to a graph realizes first the collapse of $H_{rg}^{(i)}$ to $\overline H(\G_i,\bar\e)$ and then the collapse of $\overline H(\G_i,\bar\e)$ to a graph, for each $i \in \{1,2\}$, then the vertex $v_1$ (resp. $v_2$) is an end-point of the  last (resp. first) edge in the ordering of the $4$-colored edges of $\G_1$ (resp. $\G_2$) according to Proposition \ref{CS gem-induced trisections}.

\smallskip

With regard to statement (iii),  let $\G \in G_s^{(4)}$ be a gem of the compact $4$-manifold $M^4$ with connected non-empty boundary and  let $\e$ be a cyclic permutation of $\Delta_4$, so that $ \mathcal T(\G,\e))$ is a gem-induced trisection;  since $\e_4=4$ is the unique singular color, $\partial M^4$ is represented by the 4-colored graph $\G_{\hat 4}$. Hence, Theorem \ref{thm.quasi-trisection} immediately implies: 
$$genus(\Sigma(\mathcal T(\G,\e))) = \rho_{\e_{\hat 4}}(\G_{\hat 4}) \geq \mathcal G(\partial M^4)= \mathcal H(\partial M^4).$$\qed 

 \bigskip

As a consequence of the construction described in Section \ref{sec.Kirby-diagrams} and of  Theorem \ref{th_M4(L,c)}, we are able to compute the G-trisection genus of all $\mathbb D^2$-bundles over $\mathbb S^2$ (including $\mathbb{CP}^2$  and $\mathbb S^2 \times \mathbb D^2$)  and to estimate the G-trisection genus of any $4$-manifold that is a plumbing of spheres on a linear graph. We recall that a linear plumbing of spheres can be represented by a framed link $(L_n,c)$ where $L_n$ is a chain with $n$ components (in particular, $L_2$ is the Hopf link), and $\partial M^4(L_n,c)$ - if not empty, as it happens for $M^4(L_2,(0,0))= \mathbb S^2 \times \mathbb S^2$ and $M^4(L_2,(0,1))= \mathbb S^2 \tilde \times \mathbb S^2$ - is the lens space $L(p,q)$, where $(c_1,\ldots,c_n)$ is a continued fraction expansion 
of $-\frac pq$ (see \cite{Gompf-Stipsicz}).

\medskip 

\begin{corollary}\label{linkfamilies}\ \ 
 \begin{itemize}   
  \item [(a)] Let $\xi_c$ be the $\mathbb D^2$-bundle over $\mathbb S^2$ with Euler class $c$; then\ $g_{GT}(\xi_c)=1$, for each $c \in \mathbb Z.$   
 \item [(b)] Let $M^4(L_n,c)$ be a linear plumbing of $n$ spheres ($n\in\mathbb N)$; then\ $g_{GT}(M^4(L_n,c))\leq n$, for each framing  $c\in \mathbb Z ^n.$ 
  \item [(c)] Let $\mathbb Y^4_h$ be the genus $h$ $4$-dimensional handlebody; then\ $g_{GT}(\mathbb Y^4_h)=h$, for each $h \ge 1.$  
 \end{itemize}
\end{corollary}

\dimo
Theorem \ref{th_M4(L,c)} and Corollary \ref{cor-framed-links} directly yield statement (b), as well as inequality $g_{GT}(\xi_c) \le 1$, since $\xi_c$ is the $4$-manifold associated to the $c$-framed unknot, for each $c \in \mathbb Z$; 
Proposition  \ref{g_GT properties}(i) completes the proof of statement (a).

As already pointed out in the proof of Theorem \ref{th_M4(L,c)}, for each $h\geq 1$,  the $5$-colored graph in \cite{generalized-genus} representing the genus $h$ handlebody induces a genus $h$ B-trisection for each cyclic permutation $\e$ of $\Delta_4.$ Moreover, Proposition  \ref{g_GT properties}(iii) ensures $g_{GT}(\mathbb Y^4_h) \ge h$ yielding statement (c).\qed

\begin{remark} {\em Corollary \ref{linkfamilies}(a) proves that the G-trisection genus is not finite-to-one in the boundary case, exactly as it happens for the regular genus (see \cite{generalized-genus}). Moreover, the two invariants share also the open problem about their finiteness-to-one, if a fixed (possibly empty) boundary is assumed. The problem is obviously relevant in order to understand if they are actually able to distinguish different PL structures on the same TOP $4$-manifold,  a question that is also related to that of the additivity of the G-trisection genus. 
\\
As regards this point, note that, in full analogy with regular genus, additivity of the trisection genus for closed $4$-manifolds would imply the Smooth Poincar\'e Conjecture, via a  theorem by Wall: compare  \cite[Section 6.2]{[Casali-Cristofori-Gagliardi RIMUT 2020]} and \cite{Lambert-Cole-Meier}. The clear analogies between problems and arguments   involving either of the two invariants, as well as the G-trisection genus itself, highlight their strict relation, which will be confirmed in the next section by establishing the coincidence between the G-trisection genus, the trisection genus and half the regular genus for a wide class of simply-connected closed $4$-manifolds.
}\end{remark}

\section{ Colored graphs minimizing the G-trisection genus}\label{sec.minimize_g_GT}

In order to analyze when a gem-induced trisection attains the (minimal) G-trisection genus of a compact $4$-manifold, the following lemma (similar to what already pointed out in the closed case by \cite[Section 1]{Chu-Tillman})
is very useful. 

\begin{lemma}\label{Eulercharacteristic_bis}  Let $M^4$ be a compact PL $4$-manifold with empty (resp. connected)  boundary.   For each 5-colored graph 
$\G\in G_s^{(4)}$ representing $M^4$ and for each cyclic permutation $\e$ of $\Delta_4$, the quasi-trisection $\mathcal  T(\Gamma, \e)$ of $M^4$ obtained in Theorem \ref{thm.quasi-trisection} satisfies: 
$$genus(\Sigma(\mathcal T (\G,\e)) = \bar \chi(M^4) +  g_{\e_1,\e_3,4} + g_{\e_0,\e_2,4} - \sum_{i=0}^3 g_{\hat\e_i},$$  
where  \ $\bar \chi(M^4) =  \chi(M^4) $ \ (resp. $\bar \chi(M^4) =  \chi(M^4) +1$).   
\end{lemma}

\dimo 
In \cite{generalized-genus}, the following formula for the computation of the Euler characteristic of the singular $4$-manifold $\widehat{M^4}$ is obtained:

\begin{equation}\label{Euler}\chi (\widehat{M^4})  = 2 - 2 \rho_{\varepsilon} + \sum_{i \in \Delta_4} \rho_{\varepsilon_{\hat i}}\end{equation}

\noindent where  $\e_{\hat i} = (\e_0,\ldots,\e_{i-1},\e_{i+1},\ldots,\e_4=4)$, \   $\rho_{\varepsilon}=\rho_{\varepsilon}(\Gamma)$ and $\rho_{\varepsilon_{\hat i}}=\rho_{\varepsilon_{\hat i}}(\Gamma_{\hat i}).$
                                                                                                                          
On the other hand, the same paper proves that, for each $i \in \Delta_4$,  
\begin{equation}  \label{numerospigoli(n=4)}
g_{\widehat{\e_{i-1}},\widehat{\e_{i+1}}} = g_{\e_i,\e_{i+2},\e_{i+3}} \ =  \   (g_{\widehat{\e_{i-1}}}  +  g_{\widehat{\e_{i+1}}}  -1) +     \rho_\e - \rho_{\e_{\widehat{i-1}}} - \rho_{\e_{\widehat{i+1}}}  
 \end{equation} 
 
Hence,
$$g_{\e_1,\e_3,4} + g_{\e_0,\e_2,4} = \sum_{i=0}^3 g_{\hat\e_i}  - 2 + 2\rho_\e - \sum_{i=0}^3\rho_{\e_{\hat i}}.$$  

By substituting in formula \eqref{Euler}, we have
$$\chi (\widehat{M^4}) = \rho_{\e_{\hat 4}} + \sum_{i=0}^3 g_{\hat\e_i} - g_{\e_1,\e_3,4} - g_{\e_0,\e_2,4}.  $$
 
The thesis follows by making use of claim (iv) of Theorem \ref{thm.quasi-trisection}, stating $genus(\Sigma(\mathcal T (\G,\e)) = \rho_{\e_{\hat 4}},$  together with  \ 
$\chi (\widehat{M^4}) = \bar \chi(M^4) = \begin{cases} \chi(M^4) & \text{if } \ \partial M^4 = \emptyset \\ 
                                                    \chi(M^4) +1   & \text{if } \ \partial M^4 \  \text{is connected} \end{cases}$
\qed

Generalizing \cite{[Basak-Casali 2016]}, if  \, $rk(\pi_1(M^4))= m \ge 0$ and $rk(\pi_1(\widehat M^4))= m^{\prime} \ge 0$ ($m^{\prime} \le m$), let us now introduce the symbols $t_{j,k,l}$ in the following way: 
\begin{equation}\label{g_{j,k,l}}  
g_{j,k,l}=  (g_{\hat r} + g_{\hat s} -1) + m^{\prime} + t_{j,k,l},  \ \ \text{with} \ t_{j,k,l} \ge 0,\ \{r,s\}= \Delta_4 - \{j,k,l \}, \ \ \forall j,k,l \in \Delta_3;
\end{equation} 
\begin{equation}\label{g_{j,k,4}}  
g_{j,k,4}=  (g_{\hat r} + g_{\hat s} -1) + m + t_{j,k,4},  \ \ \text{with} \ t_{j,k,4} \ge 0,\ \ \ \{r,s\}= \Delta_3 - \{j,k\}, \ \ \forall j,k \in \Delta_3.
\end{equation}  

\medskip

In case the quasi-trisection $\mathcal  T(\Gamma, \e)$  is  
actually a gem-induced trisection, the genus of the central surface is proved  to depend only on the TOP structure of $M^4$ and on the non-negative summands $ t_{\e_1,\e_3,4}$ and $t_{\e_0,\e_2,4}$:     

\begin{proposition}\label{Eulercharacteristic_ter} 
Let $M^4$ be a compact $4$-manifold with empty or connected boundary and let $\G\in G_s^{(4)}$ be a 5-colored graph representing $M^4$  such that $\mathcal T (\G,\e)$ is a gem-induced trisection ($\e$ being a cyclic permutation of $\Delta_4$). 
Then,  
 $$genus(\Sigma(\mathcal T (\G,\e)) = \beta_2(M^4) - \beta_1(M^4) + 2m +  t_{\e_1,\e_3,4} + t_{\e_0,\e_2,4}.$$
\end{proposition}

\dimo 
The proof of  Lemma \ref{Eulercharacteristic_bis} yields, for each $\G\in G_s^{(4)}$, $genus(\Sigma(\mathcal T (\G,\e)) = \chi(\widehat{M^4}) +  g_{\e_1,\e_3,4} + g_{\e_0,\e_2,4} - \sum_{i=0}^3 g_{\hat\e_i};$  
by applying equation \eqref{g_{j,k,4}},  $genus(\Sigma(\mathcal T (\G,\e)) = \chi(\widehat{M^4}) +  2m - 2 + t_{\e_1,\e_3,4} + t_{\e_0,\e_2,4}$ follows.  
Now, if  $\mathcal T(\Gamma, \e) $ is a gem-induced trisection, $\widehat{M^4}$ turns out to be simply-connected (recall Remark \ref{rem-simply-connected});  hence,  
$$genus(\Sigma(\mathcal T (\G,\e)) \ =  \ \beta_2(\widehat{M^4}) - \beta_3(\widehat{M^4}) +  2m + t_{\e_1,\e_3,4} + t_{\e_0,\e_2,4}$$   holds.

On the other hand, the exact sequence of relative homology of the pair $(\widehat{M^4},C)$, where $C$ is the cone over $\partial M^4$, together with the fact that the inclusion $(M^4,\partial M^4)\hookrightarrow (\widehat{M^4},C)$ induces an isomorphism in homology, imply 
$$H_k(\widehat{M^4})\cong H_k(\widehat{M^4},C)\cong H_k(M^4,\partial M^4)\cong H^{4-k}(M^4), \qquad \forall k\in\{2,3\}$$
\noindent where the last isomorphism comes from Poincar\'e-Lefschetz duality.

Hence $\beta_2(\widehat{M^4})=\beta_2(M^4)$ and $ \beta_3(\widehat{M^4})=\beta_1(M^4)$, from which the statement follows.   
\qed

\bigskip

 It is pointed out in \cite{Chu-Tillman} that the trisection genus of a closed $4$-manifold $M^4$ always satisfies $g_T(M^4)\geq\beta_1(M^4)+\beta_2(M^4).$ On the other hand, in \cite{SpreerTillmann} it is proved that  the trisection genus and the second Betti number actually coincide for all standard simply-connected closed $4$-manifolds. 
 
A generalization of the above inequality allows to extend this coincidence to a wider class of compact $4$-manifolds with empty or connected boundary.  In fact, as a direct consequence of Proposition \ref{Eulercharacteristic_ter}, we have:

\begin{proposition}  \label{minimal g_GT} \ \   
Let $M^4$ be a compact $4$-manifold with empty or connected boundary which admits a gem-induced trisection.
Then, $$g_{GT}(M^4)\  \ge  \beta_2(M^4)+\beta_1(M^4)+ 2\left(m-\beta_1(M^4)\right),\ \ \text{with}\ \ m = rk(\pi_1(M^4)),$$
and the equality holds if and only if a $5$-colored graph $\G\in G_s^{(4)}$ representing $M^4$ and a cyclic permutation $\e$ of $\Delta_4$  exist, so that $t_{\e_1,\e_3,4}=t_{\e_0,\e_2,4} =0$   and $\mathcal T (\G,\e)$  is a gem-induced trisection.   \\ 
\end{proposition}

Let us now take into account two particular types of crystallizations, introduced and studied in \cite{[Basak-Casali 2016]} and \cite{[Basak]} respectively\footnote{Both semi-simple and weak semi-simple crystallizations generalize the notion of {\it simple crystallizations} for closed simply-connected $4$-manifolds: see \cite{[Basak-Spreer 2016]} and \cite{[Casali-Cristofori-Gagliardi JKTR 2015]}.}; both classes have been subsequently extended to compact $4$-manifolds with empty or connected  boundary, and proved to be ``minimal" with respect to regular genus, among all graphs representing the same  $4$-manifold (see \cite{[Casali-Cristofori-Gagliardi RIMUT 2020]}). Here, we will show that they satisfy nice properties with respect to G-trisection genus, too. 

\bigskip

\begin{definition} {\em Let $M^4$ be a compact $4$-manifold, with empty or connected boundary.  A $5$-colored graph $\G\in G_s^{(4)}$  representing $M^4$   
is called {\em semi-simple}
if \, $g_{j,k,l} = 1 + m^{\prime} \ \ \forall \ j,k,l \in \Delta_3$ \, and \, $g_{j,k,4} = 1 + m \ \ \forall \ j,k \in \Delta_3,$ \, where  \, $rk(\pi_1(M^4))= m \ge 0$ and $rk(\pi_1(\widehat M^4))= m^{\prime} \ge 0$ ($m^{\prime} \le m$).  

\noindent $\G$ is called {\em weak semi-simple} with respect to a cyclic permutation $\e$ of $\Delta_4$ 
if \, $g_{\e_i,\e_{i+2},\e_{i+4}} = 1 + m \ \ \forall \ i \in \{0,2,4\}$ \, and \, $g_{\e_i,\e_{i+2},\e_{i+4}} = 1 + m^{\prime} \ \ \forall \ i \in \{1,3\}$ (where the additions in subscripts are intended in $\mathbb{Z}_{5}$).

 \smallskip
\noindent Note that, if $\G$ is weak semi-simple, then \ $g_{\hat j}=1,\ \forall \ j\in\Delta_4$, i.e. $\G$ is a crystallization of $M^4.$

\noindent 
Semi-simple (resp. weak semi-simple) crystallizations of simply-connected $4$-manifolds are often  called {\em simple} (resp. {\em weak simple}). In the following, when $m=0$, we will use the two terms indifferently.}  
\end{definition}

\medskip 

With the notations introduced in equations \eqref{g_{j,k,l}} and \eqref{g_{j,k,4}},  $\G$ turns out to be semi-simple (resp. weak semi-simple) iff $g_{\hat j}=1$ and $t_{j,k,l} = 0$  \ $\forall j,k,l \in \Delta_4$ (resp. $g_{\hat i}=1$ and $t_{\e_i,\e_{i+2},\e_{i+4}} = 0$ \ $\forall i \in \mathbb{Z}_{5}$)

\bigskip

As pointed out in \cite{generalized-genus}, by making a comparison between equations  \eqref{g_{j,k,l}} and \eqref{g_{j,k,4}} and equation  \eqref{numerospigoli(n=4)}, the following relations are proved to hold for any crystallization of 
a compact $4$-manifold with empty or connected boundary,: 

\begin{equation} \label{genus-subgenera(t)}
\begin{aligned} \rho_{\e}  - \rho_{\e_{\hat i}}  -  \rho_{\e_{\widehat {i+2}}} -  m^{\prime} = & \ t_{\e_{i-1},\e_{i+1},\e_{{i+3}}}  \ \  \ \   \forall i \in \{2,4\}  \ \ \text{and} \\
 \rho_{\e}  - \rho_{\e_{\hat i}}  -  \rho_{\e_{\widehat {i+2}}} -  m =  & \ t_{\e_{i-1},\e_{i+1},\e_{{i+3}}}   \ \  \ \   \forall i \in \{0,1,3\}. \end{aligned} \end{equation} 
  
By summing up over $i\in \Delta_4$, we obtain: 
$$ 5  \rho_{\e} - 2 \sum_{i\in \Delta_4}  \rho_{\e_{\hat i}}  - 5m  + 2( m- m^{\prime}) = \sum_{i\in \Delta_4} t_{\e_i,\e_{i+2},\e_{i+4}},$$ 
i.e. 
\begin{equation} \label{sum_subgenera(t)}
\sum_{i\in \Delta_4}  \rho_{\e_{\hat i}}  \ =  \  \frac 5 2 (\rho_{\e} - m) + ( m- m^{\prime})  - \frac 1 2  \sum_{i\in \Delta_4} t_{\e_i,\e_{i+2},\e_{i+4}}.
\end{equation}

\begin{corollary}\label{weaksemisimple}
Let $\G\in G_s^{(4)}$ be a crystallization of a compact $4$-manifold $M^4$ with empty or connected boundary,  
with \, $rk(\pi_1(M^4))= m \ge 0$ and $rk(\pi_1(\widehat M^4))= m^{\prime} \ge 0$ ($m^{\prime} \le m$).  
Then $\G$ is weak semi-simple with respect to $\e$ if and only if \  
$$ \rho_{\e_{\hat i}} = \frac 1 2 ( \rho_{\e} - m) \ \ \forall i \in  \Delta_3 \ \ \ \ \text{and}  \ \ \ \ 
\rho_{\e_{\hat 4}}  = \frac 1 2 ( \rho_{\e} - m)+  ( m- m^{\prime}).$$
\end{corollary}

\dimo 
Let us suppose that $\G$ is a weak semi-simple crystallization with respect to $\e$, then \ $t_{\e_i,\e_{i+2},\e_{i+4}}=0$ $\forall i \in \Delta_4.$
Hence, from formulas  \eqref{sum_subgenera(t)} and \eqref{genus-subgenera(t)}, \ $\sum_{i\in \Delta_4}  \rho_{\e_{\hat i}}  \ =  \  \frac 5 2 (\rho_{\e} - m) +  (m- m^{\prime})$  holds, together with \ 
$\rho_{\e}  - \rho_{\e_{\hat i}}  -  \rho_{\e_{\widehat {i+2}}} =  m^{\prime}$  $\forall i \in  \{2,4\}$ and  $\rho_{\e}  - \rho_{\e_{\hat i}}  -  \rho_{\e_{\widehat {i+2}}} =  m$  $\forall i \in  \{0,1,3\}.$

\noindent Let us now consider the five integers  $h_i$, where $h_i= \rho_{\e_{\hat i}} - \frac 1 2 ( \rho_{\e} - m)$ \ $\forall i \in  \Delta_3$ and $h_4= \rho_{\e_{\hat 4}} - \frac 1 2 ( \rho_{\e} - m) $ $-( m- m^{\prime}).$
It is easy to check that $h_i + h_{i+2} =0$ \  $\forall i\in \Delta_4$; therefore, if an integer $j\in \Delta_4$ existed, so that $h_j <0$,  then both $h_{j+2}>0$  
and $h_{j-2}>0$ would follow from the equations $\rho_{\e}  - \rho_{\e_{\hat i}}  -  \rho_{\e_{\widehat {i+2}}} =  m^{\prime}$  $\forall i \in  \{2,4\}$ and  $\rho_{\e}  - \rho_{\e_{\hat i}}  -  \rho_{\e_{\widehat {i+2}}} =  m$  $\forall i \in  \{0,1,3\}$.  
Now, $h_{j+2}>0$  and $h_{j-2}>0$ would imply both  $h_{j+4}<0$  and $h_{j-4}<0$, from the same equations. Since $\e_{j-4}=\e_{j+1}$ and $\e_{j+4}=\e_{j-1}$ are non-consecutive in $\e$, this gives a contradiction. 
As a consequence, $h_i=0$ follows, $\forall i \in \Delta_4.$
 
\medskip

\par \noindent 
On the other hand, if \  $ \rho_{\e_{\hat i}} = \frac 1 2 ( \rho_{\e} - m)$ \ $\forall i \in  \Delta_3$ and $\rho_{\e_{\hat 4}} = \frac 1 2 ( \rho_{\e} - m)+  ( m- m^{\prime})$ hold, 
then  formula  \eqref{genus-subgenera(t)} trivially implies  $t_{\e_{i-1},\e_{i+1},\e_{i+3}}=0$ $\forall i \in  \Delta_4,$ thus proving $\G$ to be weak semi-simple with respect to $\e$.\qed

\begin{proposition} \label{minimality weak semi-simple} \ \ 
Let $M^4$ be a compact $4$-manifold, with empty or connected boundary, admitting a  crystallization $\G$ which is  
weak semi-simple with respect to a cyclic permutation $\e$ of $\Delta_4$ and such that $\mathcal T(\G,\e)$ is a gem-induced trisection. \  
Then, 
$$g_{GT}(M^4)\ =\ \beta_2(M^4)+\beta_1(M^4)+ 2\left(m-\beta_1(M^4)\right) =  \frac 1 2 \, (\mathcal G(M^4)+m),$$ 
\noindent where $m = rk(\pi_1(M^4)).$
\end{proposition}

\dimo
First of all,   note that $t_{\e_1,\e_3,4}=t_{\e_0,\e_2,4} =0$  is always satisfied by crystallizations that are weak semi-simple  with respect to $\e$; hence,  if $\G\in G_s^{(4)}$ is such a  crystallization  and $\mathcal T (\G,\e)$  is a gem-induced trisection, then 
 $g_{GT}(M^4)= \beta_2(M^4)+\beta_1(M^4)+ 2\left(m-\beta_1(M^4)\right)$  in virtue of Proposition  \ref{minimal g_GT}. Moreover, Proposition \ref{trisection}, together with Corollary \ref{weaksemisimple}, enables to compute $g_{GT}(M^4)$ directly from $\G$: 
$$g_{GT}(M^4)=genus(\Sigma(\mathcal T (\G,\e)) = \rho_{\e_{\hat 4}} =  \frac 1 2 (\rho_\e +m)$$ 
Finally, it is sufficient to recall that $\rho_\e = \mathcal G(M^4)$ holds if $\G$ is a crystallization of $M^4$ which is  weak semi-simple with respect to $\e$: see \cite{[Casali-Cristofori-Gagliardi RIMUT 2020]}.
\ \qed  

 As a consequence  of  Proposition \ref{minimality weak semi-simple}, when $m=0$, we detect a large class of closed simply-connected 4-manifolds - comprehending all standard ones -, for which  the trisection genus equals the second Betti number and also coincides with half the regular genus, while Proposition \ref{minimal g_GT} , when $m=0$,  enables to detect an even larger class of such manifolds, for which the trisection genus equals the second Betti number.  Note that for both classes  the involved minimal trisections are therefore {\it efficient} according to \cite{Lambert-Cole-Meier}.

\begin{corollary} \label{minimality closed weak simple} \ \ 
$$g_T(M^4)\ =\ \beta_2(M^4)$$
for each closed  (simply-connected) $4$-manifold $M^4$  which admits a $5$-colored graph $\G\in G_s^{(4)}$ representing $M^4$ with  $t_{\e_1,\e_3,4}=t_{\e_0,\e_2,4} =0$  
($\e$ being a cyclic permutation of $\Delta_4$), so that $\mathcal T(\G,\e)$ is a gem-induced trisection. 
\par \noindent 
In particular: 
$$g_T(M^4)\ =\ \beta_2(M^4) =\frac 1 2 
 \, \mathcal G(M^4)$$      
for each closed  (simply-connected) $4$-manifold $M^4$ admitting a  crystallization $\G$  which is weak simple with respect to a cyclic permutation $\e$ of $\Delta_4$ and such that $\mathcal T(\G,\e)$ is a gem-induced trisection.
\end{corollary}

\dimo
If $M^4$ is a closed (simply-connected) 4-manifold, $g_T(M^4)\le g_{GT} (M^4)$ holds by definition. Moreover, each gem-induced trisection of $M^4$ arising from a gem $\G\in G_s^{(4)}$ with $t_{\e_1,\e_3,4}=t_{\e_0,\e_2,4} =0$  is actually a trisection where 
all 4-dimensional handlebodies are 4-balls. Hence, \cite[Lemma 6]{SpreerTillmann} proves that it is indeed a minimal genus trisection. 
Since $m=\beta_1(M^4)=0$, Proposition \ref{minimal g_GT} (resp.  Proposition \ref{minimality weak semi-simple}) 
completes the proof of the first (resp. the second) statement. 
\qed

\begin{remark}  {\em Given a compact $4$-manifold $M^4$, with empty or connected boundary, and supposing $\G\in G_s^{(4)}$ to represent $M^4$ and to satisfy $t_{\e_1,\e_3,4}=t_{\e_0,\e_2,4} =0$ for a cyclic permutation $\e$ of $\Delta_4$, if $\mathcal T(\G,\e)= (H_{0},H_{1},H_{2})$ is a gem-induced trisection, then both $H_1$ and $H_2$ are 4-balls. Hence, as we already pointed out in Remark  \ref{birman}, $\mathcal T(\G,\e)$ induces also a special Heegaard sewing of $\mathbb S^3$ or $\partial M^4$  (according to $\partial M^4$ being empty or not).\\
\noindent  If, further, $M^4$ is closed, $\mathcal T(\G,\e)$ turns out to be  efficient, i.e. it is a $(g_{T} (M^4);0)$-trisection with $g_T(M^4)\ =\ \beta_2(M^4)$. Hence edge-colored graphs are a possible useful tool to investigate the problem of the existence of infinitely many closed $4$-manifolds admitting $(g; 0)$-trisections, for some $g \ge 3$  (see  \cite[Question 4.2]{Meier}).

Note that,  closed (simply-connected) $4$-manifolds with the same second Betti number and admitting weak simple crystallizations could be a candidate set for an affirmative answer to the above question, while it is known that only a finite number of closed $4$-manifolds admitting simple crystallizations can have the same second Betti number (see \cite{[Casali-Cristofori-Gagliardi JKTR 2015]}).
}\end{remark}

\bigskip

\noindent {\bf Acknowledgements:\ }  This work was supported by GNSAGA of INDAM and by the University of Modena and Reggio Emilia, project:  {\it ``Discrete Methods in Combinatorial Geometry and Geometric Topology"}.


\begin{thebibliography}{}

\bibitem{[Ak1]} Akbulut, S.: The Dolgachev Surface  - Disproving Harer-Kas-Kirby Conjecture,
Comment. Math. Helv., 87(1), 187-241 (2012).

\bibitem{[Ak2]} Akbulut, S.: An infinite family of exotic Dolgachev surfaces without 1- and 3- handles,
Jour. of GGT, 3, 22-43 (2009).

\bibitem{Akbulut-book} Akbulut, S.: 4-manifolds, Oxford University Press 2016.

\bibitem{[Basak]}
Basak, B.: Genus-minimal crystallizations of PL 4-manifolds, Beitr. Algebra Geom., 59(1), 101-111  (2018). 

\bibitem{[Basak-Spreer 2016]}
Basak, B., Spreer, J.: Simple crystallizations of 4-manifolds,  Adv. Geom., 16(1), 111-130 (2016).

\bibitem{[Basak-Casali 2016]}
Basak,B., Casali, M.R: Lower bounds for regular genus and gem-complexity of PL 4-manifolds, Forum Mathematicum,  29 (4), 761-773 (2017). 

\bibitem{Bell-et-al} Bell, M., Hass, J., Rubinstein, J.~H., Tillmann, S.: Computing trisections of 4-manifolds, Proc. Nat. Acad. Sci. USA, 115(43), 10901-10907  (2018).

\bibitem{Birman} Birman, J.: Special Heegaard splittings of closed orientable 3-manifolds, Topology, 17, 157-166 (1978).

\bibitem{Meier_plus3} Blair, R., Cahn, P.,  Kjuchukova, A.,  Meier, J.:  A note on three-fold branched covers of $\mathbb S^4$, (2020).\ \   arXiv:1909.11788 
  


\bibitem{Casali JKTR2000} 
Casali, M.~R.: From framed links to crystallizations of bounded 4-manifolds, J. Knot Theory Ramification, 9(4), 443-458 (2000).

\bibitem{Casali_Forum2003} 
Casali, M.~R.: On the regular genus of 5-manifolds with free fundamental group,  Forum Math., 15, 465-475 (2003).  

\bibitem{generalized-genus} M.~R.~Casali - P.~Cristofori, Classifying compact 4-manifolds via generalized regular genus and G-degree, Ann. Inst. Henri Poincar\`e D (2022), 
to appear. \ \      arXiv:1912.01302

\bibitem{Casali-Cristofori Kirby-diagrams} Casali, M.~R., Cristofori, P.: Kirby diagrams and $5$-colored graphs representing compact $4$-manifolds,  (2021), submitted.  \ \ arXiv:2101.10661

\bibitem{Casali-Cristofori ElecJComb 2015}
Casali, M.~R., Cristofori, P.: Cataloguing PL 4-manifolds by gem-complexity, Electron. J. Combin., 22(4), \#P4.25 (2015).

\bibitem{Casali-Cristofori-Dartois-Grasselli}
Casali, M.R., Cristofori, P., Dartois, S., Grasselli, L.: Topology in colored tensor models  via crystallization theory,  J. Geom. Phys., 129, 142-167 (2018).  \ https://doi.org/10.1016/j.geomphys.2018.01.001 

\bibitem{Casali-Cristofori-Gagliardi Complutense 2015}
Casali,M.~R., Cristofori, P., Gagliardi, C.: Classifying PL 4-manifolds via crystallizations: results and open problems, in:  ``Mathematical Tribute to Professor Jos\'e Mar\'ia Montesinos Amilibia", Universidad Complutense Madrid (2016). \ [ISBN: 978-84-608-1684-3]

\bibitem{[Casali-Cristofori-Gagliardi JKTR 2015]}
Casali, M.R., Cristofori, P., Gagliardi, C.: PL 4-manifolds admitting simple crystallizations: framed links and regular genus, Journal of Knot Theory and its Ramifications, 25(1), 1650005 [14 pages] (2016). \  https://doi.org/10.1142/S021821651650005X

\bibitem{[Casali-Cristofori-Gagliardi RIMUT 2020]}
Casali, M.R., Cristofori, P., Gagliardi, C.:  Crystallizations of compact 4-manifolds minimizing combinatorially defined PL-invariants, Rend. Istit. Mat. Univ. Trieste 52, 431-458 (2020).  DOI: 10.13137/2464-8728/30760


\bibitem{Casali-Cristofori-Grasselli}
Casali, M.R., Cristofori, P.,  Grasselli, L.:   G-degree for singular manifolds, RACSAM,  112 (3), 693-704 (2018).  \ \ https://doi.org/10.1007/s13398-017-0456-x  

\bibitem{Casali-Gagliardi ProcAMS}
Casali, M.~R., Gagliardi, C.: Classifying PL 5-manifolds up to regular genus seven,  Proc. Amer. Math. Soc., 120,  275-283 (1994).

\bibitem{Casali-Grasselli 2019}
Casali, M.~R., Grasselli, L.: Combinatorial properties of the G-degree, Rev. Mat. Complut., 32(1), 239-254 (2019).  


\bibitem{Casella-Kuhnel} Casella, M., K\"uhnel, W.: A triangulated $K3$ surface with the minimum number of vertices, Topology, 40(4), 753-772 (2001).
 
 \bibitem{Castro-Gay-Pinzon} 
 Castro, N.~A., Gay,D.~T, Pinzon-Caicedo, J.: Trisections of 4-manifolds with boundary, Proc. Nat. Acad. Sci. USA, 
 115(43), 10861-10868 (2018). 
 
\bibitem{Castro-Ozbagci} 
Castro, N.~A., Ozbagci, B.: Trisections of 4-manifolds via Lefschetz fibrations, Math. Res. Lett., 26(2), 383-420 (2019). \ \ http://dx.doi.org/10.4310/MRL.2019.v26.n2.a3.    

\bibitem{Chu-Tillman}  
Chu, M., Tillmann, S.: Reflections on trisection genus,  Rev. Roumaine Math. Pures Appl. (Proceedings of JARCS 2017), 64(4), 395-402 (2019). 

\bibitem{Ferri-Gagliardi Proc AMS 1982}
M.~Ferri - C.~Gagliardi,  The only genus zero n-manifold is $\mathbb S^n$, Proc. Amer. Math. Soc. 85,  638-642 (1982).


\bibitem{Ferri-Gagliardi-Grasselli} 
Ferri, M., Gagliardi, C., Grasselli, L.: A graph-theoretical representation of PL-manifolds. A survey on crystallizations, Aequationes Math., 31, 121-141 (1986).
 
 \bibitem{Gagliardi 1981}
Gagliardi, C.:  Extending the concept of genus to dimension $n$, Proc. Amer. Math. Soc.,  81, 473-481 (1981).


 \bibitem{Gay-Kirby}
Gay, D., Kirby, R.: Trisecting 4-manifolds, Geom. Topol., 20, 3097-3132 (2016).

\bibitem{Gompf-Stipsicz} R. Gompf, A. Stipsicz: $4$-manifolds and Kirby calculus, Graduate Studies in Mathematics 20, Amer. Math. Soc. (1999).

\bibitem{Grasselli-Mulazzani}
Grasselli, L., Mulazzani, M.: Compact $n$-manifolds via $(n+1)$-colored graphs: a new approach, Algebra Colloq., 27(1), 95-120 (2020).   

 \bibitem{Lambert-Cole-Meier} Lambert-Cole, P. , Meier, J.: Bridge trisections in rational surfaces, Journal of Topology and Analysis, Online Ready (2020).\ \ https://doi.org/10.1142/S1793525321500047

\bibitem{[M]} 
Mandelbaum, R.: Four-dimensional topology: an introduction, Bull. Amer. Math. Soc., 2, 1-159 (1980).

\bibitem{Meier} 
Meier, J.: Trisections and spun 4-manifolds, Math. Res. Lett., 25(5), 1497-1524 (2018).  

\bibitem{Meier-Zupan} Meier, J., Zupan, A.: Bridge trisections of knotted surfaces in 4-manifolds, Proc. Nat. Acad. Sci. USA, 115(43), 10887-10893  (2018).

\bibitem{Pezzana}
Pezzana, M.: Sulla struttura topologica delle variet\`a compatte, Atti Semin. Mat. Fis. Univ. Modena, 23, 269-277 (1974).

\bibitem{Rubinstein-Tillmann} 
Rubinstein, J.~H., Tillmann, S.: Multisections of piecewise linear manifolds, Indiana Univ. Math. J.,  69(6), 2209–2239  (2020).  

\bibitem{Spreer-Kuhnel} Spreer, J., K\"uhnel, W.: Combinatorial properties of the $K3$ surface, Simplicial blowups and slicings, Exp. Math., 20(2), 201-216 (2011).

\bibitem{SpreerTillmann} 
Spreer, J., Tillmann, S.: The trisection genus of standard simply connected PL 4-manifolds, 
34nd International Symposium on Computational Geometry (SoCG 2018). In Leibniz International Proceedings in Informatics (LIPICS), vol. 99, pp. 71:1-71:13, 2018.

\end{thebibliography}
\end{document}